\documentclass{llncs}

\pdfoutput=1

\usepackage{graphicx}

\usepackage{tikz}

\usepackage{color}

\usepackage{subfigure} 

\usepackage{amsmath}   

\usepackage{amssymb}

\usepackage{float}

\usepackage{ifpdf}

\newcommand{\Cauchy}{L}
\newcommand{\Momentum}{(\mathcal L^\dagger \mathcal L)}

\newcommand{\IR}{\mathbb{R}}

\begin{document}

\frontmatter          

\title{ Efficient Gauss-Newton-Krylov momentum conservation constrained PDE-LDDMM using the band-limited vector field parameterization }

\titlerunning{Efficient Jacobi PDE-EPDiff LDDMM}     
%
\author{Monica Hernandez}

\institute{Computer Sciences Department \\ Aragon Institute on Engineering Research \\ University of Zaragoza \\ mhg@unizar.es}

\maketitle

\begin{abstract}

The class of non-rigid registration methods proposed in the framework of PDE-constrained 
Large Deformation Diffeomorphic Metric Mapping is a particularly interesting family of physically 
meaningful diffeomorphic registration methods. 
PDE-constrained LDDMM methods are formulated as constrained variational problems,
where the different physical models are imposed using the associated partial differential equations 
as hard constraints.
Inexact Newton-Krylov optimization has shown an excellent numerical accuracy and an extraordinarily 
fast convergence rate in this framework.
However, the Galerkin representation of the non-stationary velocity fields does not provide proper geodesic paths.
In a previous work, we proposed a method for PDE-constrained LDDMM parameterized in the space of initial velocity 
fields under the EPDiff equation.
The proposed method provided geodesics in the framework of PDE-constrained LDDMM, and 
it showed performance competitive to benchmark PDE-constrained LDDMM and EPDiff-LDDMM methods.
However, the major drawback of this method was the large memory load inherent to PDE-constrained LDDMM methods and 
the increased computational time with respect to the benchmark methods.
In this work we optimize the computational complexity of the method using the band-limited vector field parameterization
closing the loop with our previous works.

\end{abstract}


\section{Introduction}

Deformable image registration is the process of computing spatial transformations between different
images so that corresponding points represent the same anatomical location. 
There exists a vast literature on deformable image registration methods with differences on the transformation
characterization, regularizers, image similarity metrics, optimization methods, and additional 
constraints~\cite{Sotiras_13}.
In the last two decades, diffeomorphic registration has arisen as a powerful paradigm for deformable image
registration~\cite{Miller_04}. 
Diffeomorphisms (i.e., smooth and invertible transformations) have become fundamental inputs in Computational
Anatomy.
There exist different big families of diffeomorphic registration methods.
Our attention in the last years has been focused to PDE-constrained diffeomorphic registration due to its relevance 
in the last decade.

PDE-constrained diffeomorphic registration augments the original variational formulation with Partial 
Differential Equations (PDEs) of interest. 
The framework seems to be very appropriate for the computation of physically meaningful transformations.
The first method was proposed by Hart et al.~\cite{Hart_09}. 
In that work, the problem was formulated as a PDE-constrained control problem subject to the state PDE
and the relationship with Beg et al. LDDMM~\cite{Beg_05} was stated.
Later on, Vialard et al. proposed a PDE-constrained method parameterized on the initial 
momentum~\cite{Vialard_12}.
More recently, Mang et al. have proposed a PDE-constrained method that extended the gradient-descent 
optimization in Hart et al. approach to inexact Newton-Krylov optimization~\cite{Mang_15}.
In a previous ArXiv publication, we have bridged the gap between Vialard et al. and Mang et al. work
by proposing a novel PDE-constrained LDDMM method parameterized on the initial momentum~\cite{Hernandez_18_ArxivMICCAI}.
In addition we have faced the huge computational complexity of PDE-constrained LDDMM using the band-limited vector
field parameterization~\cite{Hernandez_18,Hernandez_18_ArxivECCV}.

This work closes the loop between our previous works~\cite{Hernandez_18,Hernandez_18_ArxivECCV,Hernandez_18_ArxivMICCAI}
by formulating the method in~\cite{Hernandez_18_ArxivMICCAI} in the space of band-limited vector fields.
This document is intended to be a self contained equation guide of all the related methods and provide the equations of 
the closing loop methods. The results section shows, as a proof of concept, the potential of our efficient method 
for Computational Anatomy applications.

\section{Related Methods}

\subsection{LDDMM}

Let $I_0$, and $I_1$ be the source and the target images defined on the image domain $\Omega \subseteq \mathbb{R}^d$.
We denote with $Diff(\Omega)$ to the Riemannian manifold of diffeomorphisms on $\Omega$.
$V$ is the tangent space of the Riemannian structure at the identity diffeomorphism, $id$.
$V$ is made of smooth vector fields on $\Omega$.
The Riemannian metric is defined from the scalar product in $V$ 
\begin{equation}
 \langle v, w \rangle_V = \langle \Cauchy v, w \rangle_{L^2} = \int_\Omega \langle \Cauchy v(x), w(x) \rangle d\Omega, 
\end{equation}
\noindent where $\Cauchy = (Id - \alpha \Delta)^s, \alpha >0, s \in \mathbb{N}$
is the invertible self-adjoint differential operator associated with the differential structure of $Diff(\Omega)$.
We denote with $K$ to the inverse of $L$.

The LDDMM variational problem is given by the minimization of the energy functional
\begin{equation}
\label{eq:LDDMM}
E(v) =  \frac{1}{2} \int_0^1 \langle \Cauchy v_t, v_t \rangle_{L^2} dt + \frac{1}{\sigma^2} \Vert I_0 \circ (\phi^{v}_{1})^{-1} - I_1\Vert_{L^2}^2. 
\end{equation}

\noindent The problem is posed in the space of time-varying smooth flows of velocity fields in $V$, ${v} \in L^2([0,1],V)$.
Given the smooth flow ${v}:[0,1] \rightarrow V$, $v_t:\Omega \rightarrow \IR^{d} \in V$, 
the diffeomorphism $\phi_1^{v}$ is defined as the solution at time $1$ to the transport equation 
$d_t \phi_t^{v} = v_t \circ \phi_t^{v}$ with initial condition $\phi_0^{v} = id$.
The transformation $(\phi^{v}_{1})^{-1}$ computed from the minimum of $E({v})$ is the diffeomorphism that 
solves the LDDMM registration problem between $I_0$ and $I_1$.
The optimization of Equation~\ref{eq:LDDMM} was originally approached in~\cite{Beg_05} using gradient-descent in $L^2([0,1],V)$, 
yielding the update equation 

\begin{equation}
v_t^{n+1} = v_t^n - \epsilon (\nabla_v E(v))_t. 
\end{equation}

\subsection{EPDiff-LDDMM}

The geodesics of $Diff(\Omega)$ under the right-invariant Riemannian metric are uniquely determined by 
the time-varying flows of velocity fields that satisfy the Euler-Poincar\'e equation (EPDiff)~\cite{Holm_98}
\begin{equation}
\label{eq:EPDiff}
\partial_t v = -ad_v^\dagger v = - K ad_v^* L v = - K [ (Dv)^T Lv + D(Lv) v + L v \nabla \cdot v ].
\end{equation}
\noindent with initial condition $v_0 \in V$.

LDDMM can be posed in the space of initial velocity fields
 
\begin{equation}
\label{eq:EPDiffLDDMM}
E(v_0) = \frac{1}{2} \langle \Cauchy v_0, v_0 \rangle_{L^2} + \frac{1}{\sigma^2} \Vert I_0 \circ (\phi^{v}_{1})^{-1} - I_1 \Vert_{L^2}^2, 
\end{equation}

\noindent where $(\phi^{v}_{1})^{-1}$ is the solution at time $1$ to the transport equation of the flow $v_t$ that satisfies 
the EPDiff equation for $v_0$. 
The optimization of Equation~\ref{eq:EPDiffLDDMM} was originally approached using gradient-descent in $V$~\cite{Younes_07}
\begin{equation}
v_0^{n+1} = v_0^n - \epsilon \nabla_{v_0} E(v_0). 
\end{equation}

More recently, it has been proposed in~\cite{Zhang_15} to compute the gradient at $t = 1$ and to integrate backward
the reduced adjoint Jacobi field equations~\cite{Bullo_95} 
\begin{eqnarray}
\label{eq:AdjointJacobiField}
\partial_t U_t + ad_{v_t}^\dagger U_t = 0 \textnormal{ in } \Omega \times [0,1) \\
\partial_t \delta v_t + U_t - ad_{v_t} \delta v_t + ad_{\delta v_t}^\dagger v_t = 0 \textnormal{ in } \Omega \times [0,1)
\end{eqnarray}
\noindent with initial conditions $U(1) = \nabla_{v_1} E(v_0)$ and $\delta v(1) = 0$,
to get the gradient update at $t = 0$, 

\begin{equation}
v_0^{n+1} = v_0^n - \epsilon \delta v(0).  
\end{equation}

\subsection{PDE-LDDMM subject to the state equation}

The PDE-constrained LDDMM variational problem is given by the minimization of 
\begin{equation}
\label{eq:LS-LDDMM}
E(v) = \frac{1}{2} \int_0^1 \langle \Cauchy v_t, v_t \rangle_{L^2} dt + \frac{1}{\sigma^2} \Vert m(1) - I_1\Vert_{L^2}^2, 
\end{equation}
\noindent subject to the state equation
\begin{equation}
\label{eq:StateEquation}
\partial_t m(t) + \nabla m(t) \cdot v_t = 0 \textnormal{ in } \Omega \times (0,1],
\end{equation}
\noindent with initial condition $m(0) = I_0$. 
The compressible PDE-constrained problem was proposed by Hart et al. with gradient-descent 
optimization~\cite{Hart_09}.
Mag et al. introduced the incompressibility constraint and solved the problem using inexact 
Newton-Krylov optimization~\cite{Mang_15}.
 
In PDE- LDDMM, the gradient and the Hessian are computed using the method of Lagrange multipliers. 
Thus, we define the Lagrange multiplier $\lambda:\Omega \times [0,1] \rightarrow \mathbb{R}$ associated with 
the state equation, and we build the augmented Lagrangian
\begin{equation}
\label{eq:LS-LDDMM-Aug}
E_{aug}(v) =  E(v) +  
\int_0^1 \langle \lambda(t), \partial_t m(t) + \nabla m(t) \cdot v_t \rangle_{L^2} dt. \nonumber
\end{equation}

\noindent The first-order variation of the augmented Lagrangian yields the expression of the gradient 
\begin{eqnarray}
\partial_t m(t) + \nabla m(t) \cdot v_t = 0 \textnormal{ in } \Omega \times (0,1] \label{eq:FirstOrder1} \\
-\partial_t \lambda(t) - \nabla \cdot ( \lambda(t) \cdot v_t ) = 0 \textnormal{ in } \Omega \times [0,1) \label{eq:FirstOrder2}  \\
(\nabla_v E_{aug}(v))_t = \Momentum v_t + \lambda(t) \cdot \nabla m(t) \textnormal{ in } \Omega \times [0,1] \label{eq:FirstOrder4}  
\end{eqnarray}

\noindent subject to the initial and final conditions 
$m(0) = I_0$ and $\lambda(1) = -\frac{2}{\sigma^2}(m(1) - I_1) \textnormal{ in } \Omega$. 
In the following, we will recall $m$ as the state variable and $\lambda$ as the adjoint variable.
Equations~\ref{eq:FirstOrder1} and~\ref{eq:FirstOrder2} will be recalled as the state and adjoint equations, respectively.

The second-order variation of the augmented Lagrangian yields the expression of the Hessian-vector product
\begin{multline}
(H_v E_{aug}(v))_t \cdot \delta v(t) = \\
\Momentum \delta v(t) + \delta \lambda(t) \cdot \nabla m(t) - \lambda(t) \cdot \nabla \delta m(t) \textnormal{ in } \Omega \times [0,1] \label{eq:SecondOrder3}   
\end{multline}
\noindent where
\begin{eqnarray}
\label{eq:SecondOrderVariation}
\partial_t \delta m (t) + \nabla \delta m(t) \cdot v_t + \nabla m(t) \cdot \delta v(t) = 0 \textnormal{ in } \Omega \times (0,1] \label{eq:SecondOrder1}\\
-\partial_t \delta \lambda(t) - \nabla \cdot ( \delta \lambda(t) \cdot v_t ) + \nabla \cdot ( \lambda(t) \cdot \delta v(t) ) = 0 \textnormal{ in } \Omega \times [0,1) \label{eq:SecondOrder2}
\end{eqnarray}
\noindent subject to 
$\delta m(0) = 0$ and $\delta \lambda(1) = -\frac{2}{\sigma^2}\delta m(1) \textnormal{ in } \Omega$.
Equation~\ref{eq:SecondOrder1} corresponds with the incremental state equation.
Equation~\ref{eq:SecondOrder2} corresponds with the incremental adjoint equation.

\subsection{PDE-LDDMM subject to the deformation state equation}

This variant of PDE-constrained LDDMM is formulated from the minimization of Equation~\ref{eq:LS-LDDMM}
subject to the deformation state equation
\begin{equation}
\partial_t \phi(t) + D \phi(t) \cdot v_t = 0 \textnormal{ in } \Omega \times (0,1],
\end{equation}
\noindent and the incompressibility constraint
\begin{equation}
\gamma \nabla \cdot v_t = 0 \textnormal{ in } \Omega \times [0, 1].
\end{equation} 
\noindent The compressible PDE-constrained problem was proposed by Polzin et al. with gradient-descent 
optimization~\cite{Polzin_16}.

The Lagrange multipliers are $\rho:\Omega \times [0,1] \rightarrow \mathbb{R}^d$, associated with 
the deformation state equation, and $p: \Omega \times [0,1] \rightarrow  \mathbb{R}^d$, associated with the incompressibility 
constraint.
The augmented Lagrangian is given by
\begin{equation}
\label{eq:LS-LDDMM-Aug-Polzin}
E_{aug}(v) =  E(v) +  
\int_0^1 \langle \rho(t), \partial_t \phi(t) + D \phi(t) \cdot v_t \rangle_{L^2} dt. \nonumber
\end{equation}

The expression of the gradient is given by the first-order variation of the augmented Lagrangian
\begin{eqnarray}
\partial_t \phi(t) + D \phi(t) \cdot v_t = 0 \textnormal{ in } \Omega \times (0,1] \label{eq:PFirstOrder1} \\
-\partial_t \rho(t) - \nabla \cdot (\rho(t) \cdot v_t ) = 0 \textnormal{ in } \Omega \times [0,1) \label{eq:PFirstOrder2} \\
(\nabla_v E_{aug}(v))_t = \Momentum v_t + D \phi(t) \cdot \rho(t) \textnormal{ in } \Omega \times [0,1] \label{eq:PFirstOrder3}
\end{eqnarray}
\noindent subject to the initial and final conditions 
$\phi(0) = id$, and $\rho(1) = \lambda(1) \cdot \nabla m(1)$. 

From the second-order variation of the augmented Lagrangian, we obtain the expression of the Hessian-vector product 
\begin{multline}
(H_v E_{aug}(v))_t \cdot \delta v(t) = \\
\Momentum \delta v(t) + D \delta \phi(t) \cdot \rho(t) -  D \phi(t) \cdot \delta \rho(t) \textnormal{ in } \Omega \times [0,1] \label{eq:PSecondOrder3} 
\end{multline}

\noindent where

\begin{eqnarray}
\partial_t \delta \phi(t) + D \delta \phi(t) \cdot v_t + D \phi(t) \cdot \delta v(t) = 0 \textnormal{ in } \Omega \times (0,1] \label{eq:PSecondOrder1} \\
-\partial_t \delta \rho(t) - \nabla \cdot ( \delta \rho(t) \cdot v_t ) + \nabla \cdot (\rho(t) \cdot \delta v(t)) = 0 \textnormal{ in } \Omega \times [0,1) \label{eq:PSecondOrder2} 
\end{eqnarray}

\noindent subject to $\delta \phi(0) = 0$, $\delta \rho(1) = \delta \lambda(1) \cdot \nabla m(1) - \lambda(1) \cdot \delta \nabla m(1)$.

\subsection{Jacobi PDE-EPDiff LDDMM subject to the state equation}

In~\cite{Hernandez_18_ArxivMICCAI} we proposed bridging the gap between momentum conservation constrained LDDMM and 
the PDE-constrained LDDMM method in~\cite{Mang_15} using the adjoint Jacobi field equations.
With this approach, the integration of the adjoint equation is not needed.
During the derivation of the equations we first explored the idea of transporting the gradient and the Hessian-vector
product using the adjoint Jacobi equations.
However, we found that the resulting method did not converge.
We found that transporting the vectors differently (i.e. using the adjoint Jacobi equations for the gradient, 
and the incremental adjoint Jacobi equations for the Hessian-vector products) yields the desired convergence behavior.

Thus, the PDE-constrained problem is given by the minimization of the energy functional
\begin{equation}
 E(v_0) = \frac{1}{2} \langle L v_0, v_0 \rangle_{L^2} + \frac{1}{\sigma^2} \Vert m(1) - I_1 \Vert_{L^2}^2,
\end{equation}
\noindent subject to the EPDiff and the state equations
\begin{eqnarray}
\partial_t v_t + ad_{v_t}^\dagger v_t = 0 \textnormal{ in } \Omega \times (0,1] \\
\partial_t m(t) + \nabla m(t) \cdot v_t = 0 \textnormal{ in } \Omega \times (0,1],
\end{eqnarray}

\noindent with initial conditions $v(0) = v_0$ and $m(0) = I_0$, respectively.

Optimization is performed combining the method of Lagrange multipliers with inexact Gauss-Newton-Krylov methods
in the following way. Let
$w: \Omega \times [0,1] \rightarrow \mathbb{R}^d$ and
$\lambda: \Omega \times [0,1] \rightarrow \mathbb{R}$ 
be the Lagrange multipliers associated with the EPDiff and the state equations.
We build the augmented Lagrangian
\begin{multline}
\label{eq:MCC-LDDMM-Aug}
E_{aug}(v_0) =  E(v_0) + \int_0^1 \langle w(t), \partial_t v(t) + ad_{v_t}^\dagger v_t \rangle_{L^2} dt \\
+ \int_0^1 \langle \lambda(t), \partial_t m(t) + \nabla m(t) \cdot v_t \rangle_{L^2} dt.
\end{multline}

Similarly to~\cite{Zhang_15}, the gradient is computed at $t = 1$, $\nabla_{v_1} E(v_0) = \lambda(1) \cdot \nabla m(1)$
and integrated backward using the reduced adjoint Jacobi field equations (Equation~\ref{eq:AdjointJacobiField})
to obtain $\nabla_{v_0} E(v_0)$.

The second-order variations of the augmented Lagrangian on $w$ and $\lambda$ yield the incremental EPDiff and 
incremental state equations, needed for the computation of the Hessian-vector product. Thus,
\begin{eqnarray}
 \partial_t \delta v_t + ad_{\delta {v_t}}^\dagger v_t + ad_{v_t}^\dagger \delta v_t= 0 \textnormal{ in } \Omega \times (0,1] \\
 \partial_t \delta m(t) + \nabla \delta m(t) \cdot v_t + \nabla m(t) \cdot \delta v_t = 0 \textnormal{ in } \Omega \times (0,1]
\end{eqnarray}
\noindent with initial conditions $\delta v(0) = 0$ and $\delta m(0) = 0$.

The Hessian-vector product $H_{v_0} E(v_0) \cdot \delta {v_0}$ is computed from the Hessian-vector product at $t = 1$, which is integrated backward
using the reduced incremental adjoint Jacobi field equations
\begin{eqnarray}
 \partial_t \delta U + ad_{\delta v}^\dagger U + ad_v^\dagger \delta U= 0 \textnormal{ in } \Omega \times [0,1) \\
 \partial_t \delta w + \delta U - ad_{\delta v} w - ad_v \delta w + ad_{\delta w}^\dagger v + ad_{w}^\dagger \delta v = 0 \textnormal{ in } \Omega \times [0,1)
\end{eqnarray}

\noindent with initial conditions $\delta U(1) = K( \delta \lambda(1) \cdot \nabla m(1) - \lambda(1) \cdot \nabla \delta m(1) )$, 
and $\delta w(1) = 0$. 

\subsection{Jacobi PDE-EPDiff LDDMM subject the deformation state equation}

In~\cite{Hernandez_18_ArxivECCV} we explored the behavior of different variants of the PDE-constrained LDDMM problem with the band-limited vector field parameterization.
The best performing method was the Newton-Krylov extension of the method proposed in~\cite{Polzin_14}.
In this work we provide the equations of the Jacobi PDE-EPDiff LDDMM version of the method.

The PDE-constrained problem is given by the minimization of the energy functional
\begin{equation}
 E(v_0) = \frac{1}{2} \langle L v_0, v_0 \rangle_{L^2} + \frac{1}{\sigma^2} \Vert m(1) - I_1 \Vert_{L^2}^2,
\end{equation}
\noindent subject to the EPDiff and the deformation state equations
\begin{eqnarray}
\partial_t v_t + ad_{v_t}^\dagger v_t = 0 \textnormal{ in } \Omega \times (0,1] \\
\partial_t \phi(t) + D \phi(t) \cdot v_t = 0 \textnormal{ in } \Omega \times (0,1],
\end{eqnarray}

\noindent with initial conditions $v(0) = v_0$ and $\phi(0) = id$, respectively.
The state variable $m$ is computed from $m(t) = I_0 \circ \phi(t)$.

The augmented Lagrangian is given by
\begin{multline}
\label{eq:MCC-LDDMM-Aug}
E_{aug}(v_0) =  E(v_0) + \int_0^1 \langle w(t), \partial_t v(t) + ad_{v_t}^\dagger v_t \rangle_{L^2} dt \\
+ \int_0^1 \langle \rho(t), \partial_t \phi(t) + D \phi(t) \cdot v_t \rangle_{L^2} dt.
\end{multline}

The gradient is computed at $t = 1$, $\nabla_{v_1} E(v_0) = D \phi(1) \cdot \rho(1)$
and integrated backward using the reduced adjoint Jacobi field equations (Equation~\ref{eq:AdjointJacobiField})
to obtain $\nabla_{v_0} E(v_0)$.
Thus, 

\begin{eqnarray}
\partial_t \phi(t) + D \phi(t) \cdot v_t = 0 \textnormal{ in } \Omega \times (0,1], \\
\partial_t U_t + ad_{v_t}^\dagger U_t = 0 \textnormal{ in } \Omega \times [0,1) \\
\partial_t \delta v_t + U_t - ad_{v_t} \delta v_t + ad_{\delta v_t}^\dagger v_t = 0 \textnormal{ in } \Omega \times [0,1)
\end{eqnarray}

\noindent where $U(1) = K( D \phi(1) \cdot \rho(1) )$.

The Hessian-vector product $H_{v_0} E( v_0 ) \cdot \delta {v_0}$ is computed from the Hessian-vector product at $t = 1$, which is integrated backward
using the reduced incremental adjoint Jacobi field equations

\begin{eqnarray}
 \partial_t \delta \phi(t) + D \delta \phi(t) \cdot v_t + D \phi(t) \cdot \delta v_t = 0 \textnormal{ in } \Omega \times (0,1], \\
 \partial_t \delta U_t + ad_{\delta v_t}^\dagger U_t + ad_{v_t}^\dagger \delta U_t = 0 \textnormal{ in } \Omega \times [0,1) \\
 \partial_t \delta w_t + \delta U_t - ad_{\delta v_t} w_t - ad_{v_t} \delta w_t + ad_{\delta w_t}^\dagger v_t + ad_{w_t}^\dagger \delta v_t = 0 \textnormal{ in } \Omega \times [0,1)
\end{eqnarray}

\noindent with initial conditions $\delta U(1) = K( D \delta \phi(1) \cdot \rho(1) - D \phi(1) \cdot \delta \rho(1) ) $, and $\delta w(1) = 0$.

\subsection{Gauss-Newton-Krylov optimization}

By construction, the Hessian is positive definite in the proximity of a local minimum.
However, it can be indefinite or singular far away from the solution.
In this case, the search directions obtained with PCG are not guaranteed to be descent directions. 
In order to overcome this problem, one can use a Gauss-Newton approximation dropping expressions of $H_v E(v_0) \cdot \delta  v_0$
to guarantee that the matrix is definite positive.

The minimization using a second-order inexact Gauss-Newton-Krylov method yields to the update equation
\begin{equation}
\label{eq:Newton}
 v_0^{n+1} = v_0^n - \epsilon \delta v_0^n,
\end{equation}
\noindent where $ \delta v_0^n$ is computed from PCG on the system 
\begin{equation}
\label{eq:KKT}
H_{v_0} E( v_0^n) \cdot \delta {v_0}^n = - \nabla_{v_0} E(v_0^n).
\end{equation}
\noindent In this work, we consider CG with the gradient and the Hessian computed on $V$ instead of $L^2$.

\section{Methods parameterized in the space of band-limited vector fields}

\subsection{Background on the space of band-limited vector fields}

Let $\widetilde{\Omega}$ be the discrete Fourier domain truncated with frequency bounds $K_1,$ $\dots,$ $K_d$.
We denote with $\widetilde{V}$ the space of discretized band-limited vector fields on $\Omega$ with these frequency bounds. 
The elements in $\widetilde{V}$ are represented in the Fourier domain as $\tilde{v}: \widetilde{\Omega} \rightarrow \mathbb{C}^d$,
$\tilde{v}(k_1, \dots, k_d)$, and in the spatial domain as $\iota(\tilde{v}):\Omega \rightarrow \mathbb{R}^d$,
\begin{equation}
\iota(\tilde{v})(x_1, \dots, x_d) = \sum_{k_1 = 0}^{K_1} \dots \sum_{k_d = 0}^{K_d} \tilde{v}(k_1, \dots, k_d) e^{2 \pi i k_1 x_1} \dots e^{2 \pi i k_d x_d}. 
\end{equation}

\noindent The application $\iota:\widetilde{V} \rightarrow V$ denotes the natural inclusion mapping of $\widetilde{V}$ in $V$.
The aplication $\pi: V \rightarrow \widetilde{V}$ denotes the projection of $V$ onto $\widetilde{V}$.

The space of band-limited vector fields $\widetilde{V}$ has a finite-dimensional Lie algebra structure using the truncated 
convolution in the definition of the Lie bracket~\cite{Zhang_15}.
We denote with $Diff(\widetilde{\Omega})$ to the finite-dimensional Riemannian manifold of diffeomorphisms on $\widetilde{\Omega}$ 
with corresponding Lie algebra $\widetilde{V}$.
The Riemannian metric in $Diff(\widetilde{\Omega})$ is defined from the scalar product
\begin{equation}
 \langle \tilde{v}, \tilde{w} \rangle_{\tilde{V}} = \langle \tilde{\Cauchy} \tilde{v}, \tilde{w} \rangle_{l^2},  
\end{equation}
\noindent where $\tilde{\Cauchy}$ is the projection of operator $\Cauchy$ in the truncated Fourier domain. 
Similarly, we will denote with $\tilde{K}$, $\widetilde{\nabla}$, and $\widetilde{\nabla \cdot}$ to 
the projection of operators $K$, $\nabla$, and $\nabla \cdot$ in the truncated Fourier domain.
In addition, we will denote with $\star$ to the truncated convolution.

The EPDiff-equation in the space of band-limited vector fields is given by
\begin{equation}
 \partial_t \tilde{v}_t = \widetilde{ad}^\dagger_{\tilde{v}_t} \tilde{v}_t = -\tilde{K} [ (\tilde{D}\tilde{v})^T \star \tilde{L}\tilde{v} + \tilde{D}(\tilde{L}\tilde{v}) \star \tilde{v} + \tilde{L} \tilde{v} \nabla \star \tilde{v} ]. 
\end{equation}
The adjoint operator is given by
\begin{equation}
 \widetilde{ad}_{\tilde{v}}{\tilde{w}} = \tilde{D}\tilde{v} \star \tilde{w} - \tilde{D}\tilde{w} \star \tilde{v}.
\end{equation}

\subsection{BL Jacobi PDE-EPDiff LDDMM subject to the state equation}

The variational problem is given by the minimization of 
\begin{equation}
\label{eq:BLEnergy}
 E(\tilde{v}_0) = \frac{1}{2} \langle \tilde{L}\tilde{v}_0, \tilde{v}_0 \rangle_{l^2} + \frac{1}{\sigma^2} \Vert m(1) - I_1 \Vert_{L^2}^2
\end{equation}
\noindent subject to 
\begin{eqnarray}
\partial_t \tilde{v}_t + \widetilde{ad}^\dagger_{\tilde{v}_t} \tilde{v}_t = 0 \textnormal{ in } \Omega \times (0,1] \\
\partial_t m(t) + \nabla m(t) \cdot \iota(\tilde{v}_t) = 0 \textnormal{ in } \Omega \times (0,1],
\end{eqnarray}
\noindent with initial conditions $\tilde{v}(0) = \tilde{v}_0$ and $m(0) = I_0$. 

The expression of the gradient is computed from the reduced adjoint Jacobi field equations in the space of band-limited vectors yielding
\begin{eqnarray}
\partial_t \tilde{U}_t + \widetilde{ad}^\dagger_{\tilde{v}_t} \tilde{U}_t = 0 \textnormal{ in } \Omega \times [0,1) \\
\partial_t \delta\tilde{v}_t + \tilde{U} - \widetilde{ad}^\dagger_{\tilde{v}_t} \delta \tilde{v}_t + \widetilde{ad}^\dagger_{\delta \tilde{v}_t} \tilde{v}_t = 0 \textnormal{ in } \Omega \times [0,1) 
\end{eqnarray}
\noindent where  
\begin{eqnarray}
\tilde{U}(1) = \tilde{K}( \pi( \lambda(1) \cdot \nabla m(1) ) ).
\end{eqnarray}

The expression of the Hessian-vector product is computed from the reduced incremental adjoint Jacobi equations in the space of band-limited vector fields
\begin{eqnarray}
\partial_t \delta \tilde{U}_t + \widetilde{ad}^\dagger_{\tilde{\delta v}_t} \tilde{U}_t + \widetilde{ad}^\dagger_{\tilde{v}_t} \delta \tilde{U}_t = 0 \textnormal{ in } \Omega \times [0,1) \\
\partial_t \delta \tilde{w}_t + \delta \tilde{U}_t - \widetilde{ad}_{\delta \tilde{v}_t} \tilde{w}_t - \widetilde{ad}_{\tilde{v}_t} \delta \tilde{w}_t + \widetilde{ad}_{\delta \tilde{w}_t}^\dagger \tilde{v}_t + \widetilde{ad}_{\tilde{w}_t}^\dagger \delta \tilde{v}_t = 0 \textnormal{ in } \Omega \times [0,1)
\end{eqnarray}
\noindent where  
\begin{eqnarray}
\delta \tilde{U}(1) = \pi( \tilde{K}(\delta \lambda(1) \cdot \nabla m(1) - \lambda(1) \cdot \nabla \delta m(1)) ).
\end{eqnarray}

\subsection{BL Jacobi PDE-EPDiff LDDMM subject to the deformation state equation}

The variational problem is given by the minimization of Equation~\ref{eq:BLEnergy} subject to 
\begin{eqnarray}
\partial_t \tilde{v}_t + \widetilde{ad}^\dagger_{\tilde{v}_t} \tilde{v}_t = 0 \textnormal{ in } \Omega \times [0,1) \\
\partial_t \tilde{\phi}(t) + \widetilde{D} \tilde{\phi}(t) \star \tilde{v}_t = 0 \textnormal{ in } \Omega \times [0,1).
\end{eqnarray}
The expression of the gradient is given by
\begin{eqnarray}
\partial_t \tilde{v}_t + \widetilde{ad}^\dagger_{\tilde{v}_t} \tilde{v}_t = 0 \textnormal{ in } \Omega \times [0,1) \\
\partial_t \tilde{\phi}(t) + \widetilde{D} \tilde{\phi}(t) \star \tilde{v}_t = 0 \textnormal{ in } \Omega \times (0,1]  \\
\partial_t \tilde{U}_t + \widetilde{ad}^\dagger_{\tilde{v}_t} \tilde{U}_t = 0 \textnormal{ in } \Omega \times (0,1] \\
\partial_t \delta \tilde{v} + \tilde{U} - \widetilde{ad}_{\tilde{v}} \delta \tilde{v} + \widetilde{ad}^\dagger_{\delta \tilde{v}} \tilde{v} = 0 \textnormal{ in } \Omega \times [0,1)
\end{eqnarray}

\noindent with initial conditions $\tilde{U}(1) = \tilde{K}(\tilde{D} \tilde{\phi}(1) \star \tilde{\rho}(1))$.

On the other hand, the expression of the Hessian-vector product is given by
\begin{eqnarray}
\partial_t \delta \tilde{v}_t + \widetilde{ad}^\dagger_{\tilde{\delta v}_t} \tilde{v}_t + \widetilde{ad}^\dagger_{\tilde{v}_t} \tilde{\delta v}_t = 0 \textnormal{ in } \Omega \times [0,1) \\
\partial_t \delta \tilde{\phi}(t) + \widetilde{D} \delta \tilde{\phi}(t) \star \tilde{v}_t + \widetilde{D} \tilde{\phi}(t) \star \delta \tilde{v}(t) = 0  \\
\label{eq:IncrementalDeformationAdjoint}
\partial_t \delta \tilde{U}_t + \widetilde{ad}^\dagger_{\tilde{\delta v}_t} \tilde{U}_t + \widetilde{ad}^\dagger_{\tilde{v}_t} \delta \tilde{U}_t = 0 \textnormal{ in } \Omega \times [0,1) \\
\partial_t \delta \tilde{w}_t + \delta \tilde{U}_t - \widetilde{ad}_{\delta \tilde{v}_t} \tilde{w}_t - \widetilde{ad}_{\tilde{v}_t} \delta \tilde{w}_t + \widetilde{ad}_{\delta \tilde{w}_t}^\dagger \tilde{v}_t + \widetilde{ad}_{\tilde{w}_t}^\dagger \delta \tilde{v}_t = 0 \textnormal{ in } \Omega \times [0,1)
\end{eqnarray}
\noindent where  
\begin{equation}
\delta \tilde{U}(1) = \tilde{K}(\tilde{D} \delta \tilde{\phi}(1) \star \tilde{\rho}(1) - \tilde{D} \tilde{\phi}(1) \star \delta \tilde{\rho}(1)).
\end{equation}

\section{Results}

\begin{figure} [!t]
\centering
\scriptsize
\begin{tabular}{ccc}
\includegraphics[width = 0.3 \textwidth]{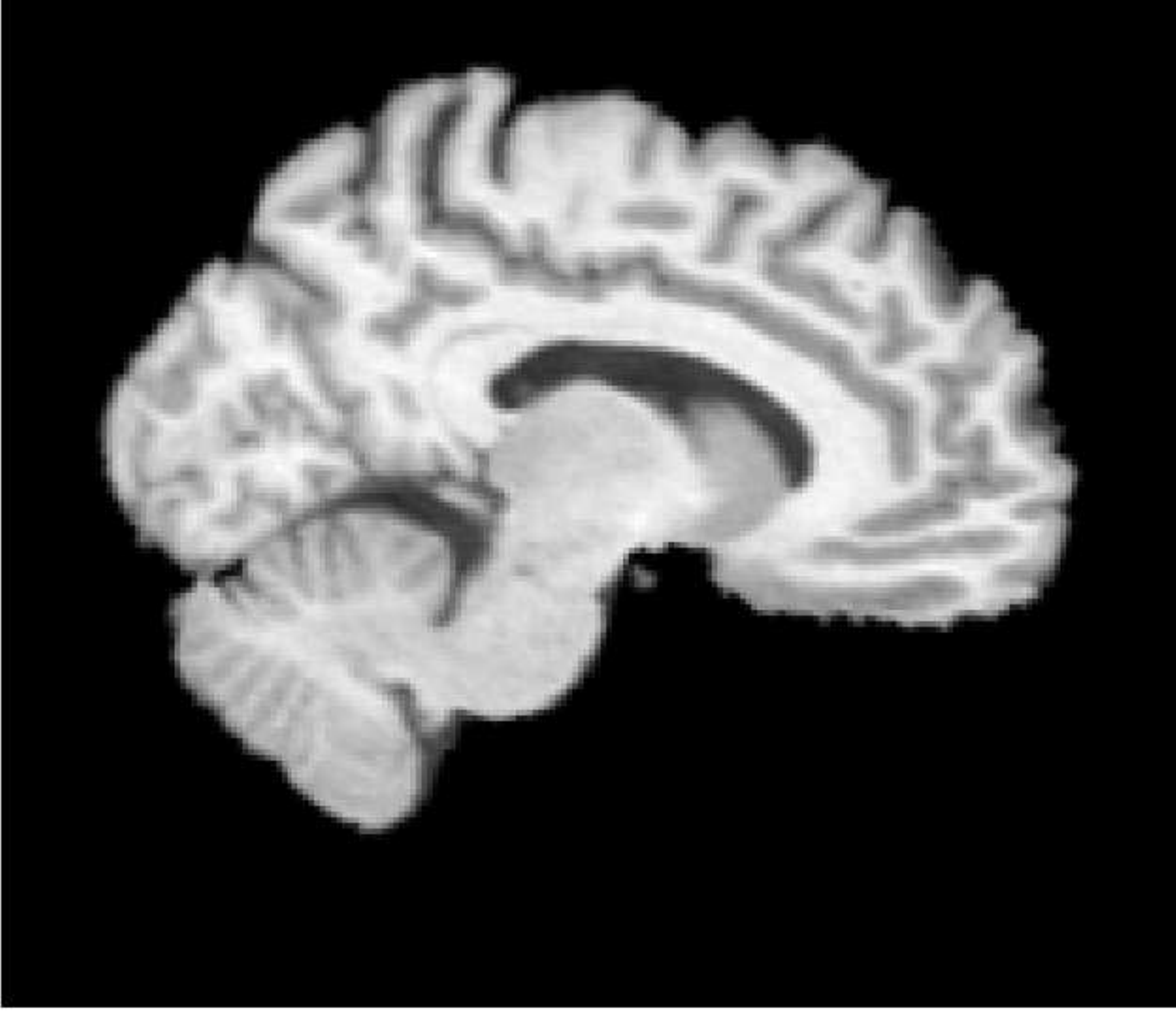} 
&
\includegraphics[width = 0.3 \textwidth]{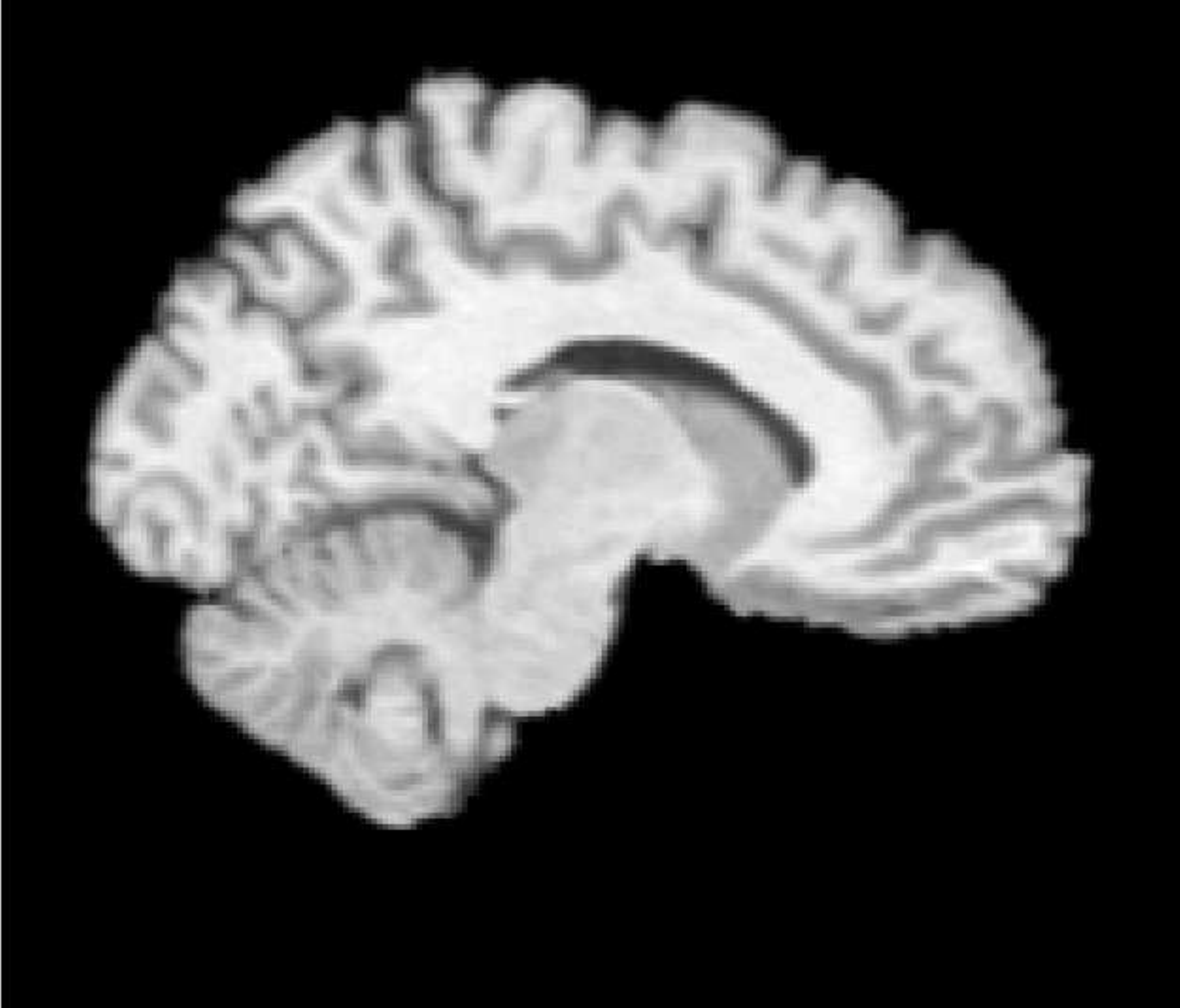} 
&
\includegraphics[width = 0.3 \textwidth]{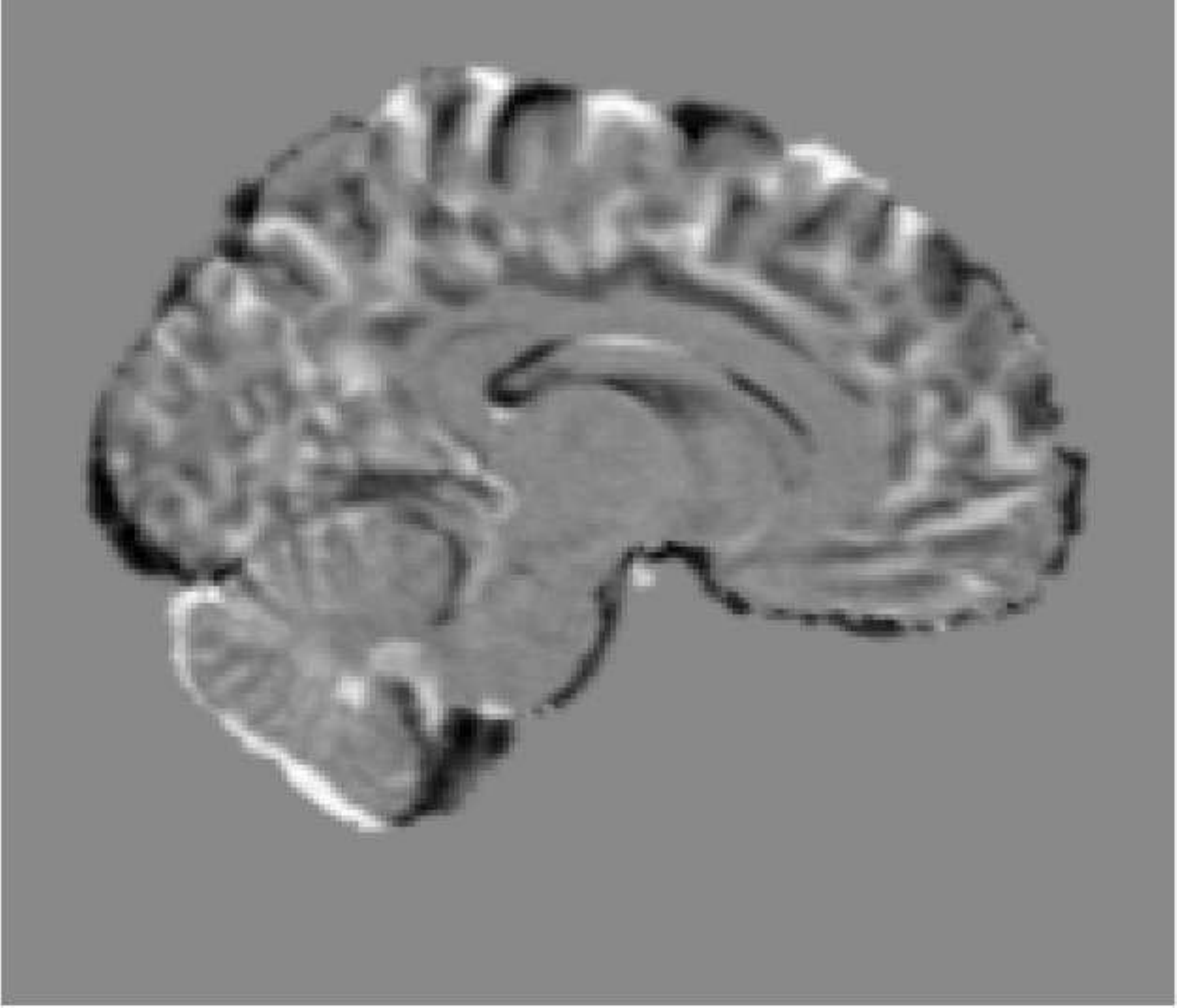} 
\end{tabular}
\caption{\small Source, target and difference before registration.}
\label{fig:sources}
\end{figure}

The experiments have been conducted on the Non-rigid Image Registration Evaluation Project database (NIREP) with volumes of size $180 \times 210 \times 180$. 
Figure~\ref{fig:sources} shows the source and target images and the differences before registration.

\begin{figure} [!t]
\centering
\scriptsize
\begin{tabular}{cccc}
\includegraphics[width = 3 cm, height = 2.5 cm]{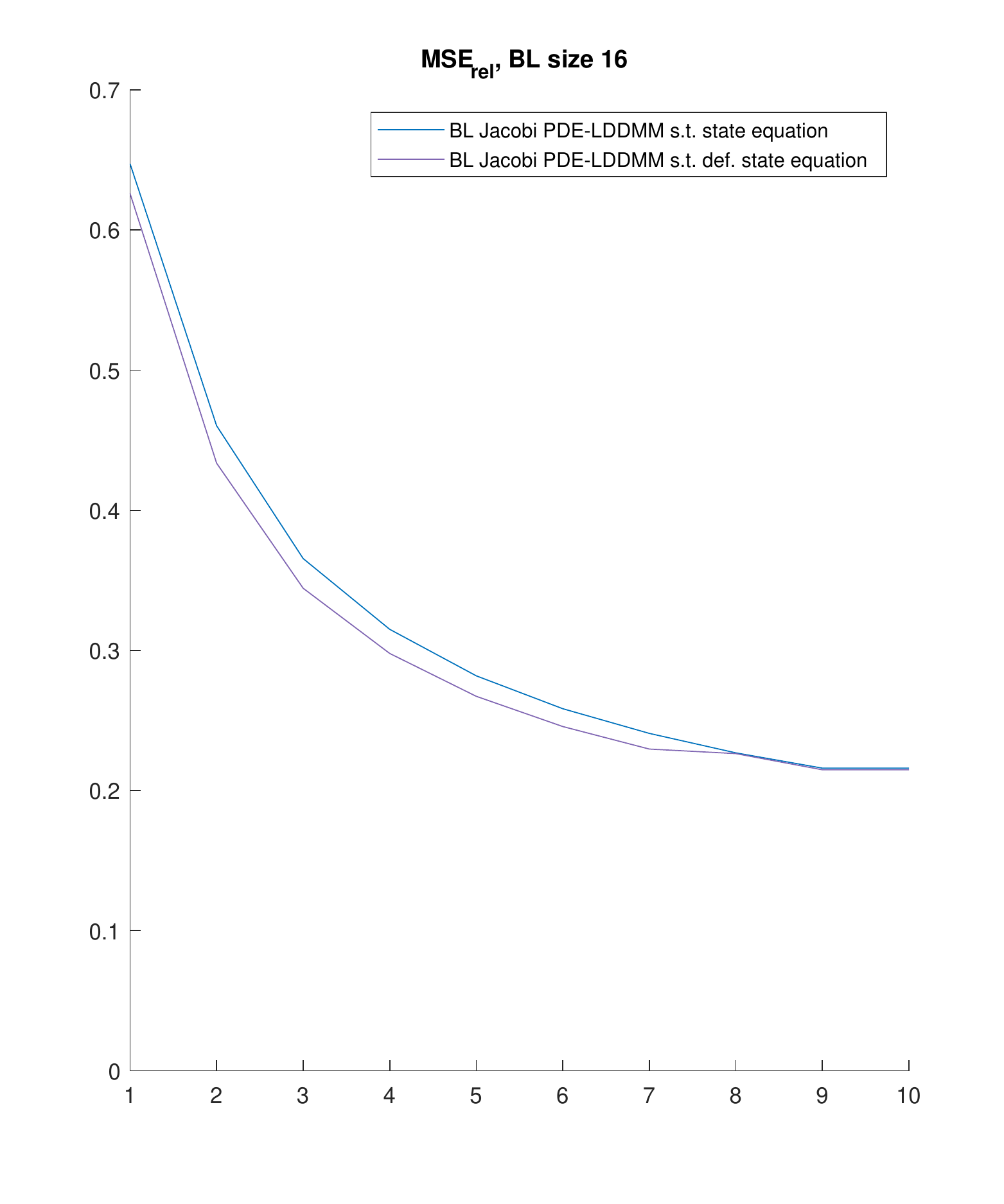} 
&
\includegraphics[width = 3 cm, height = 2.5 cm]{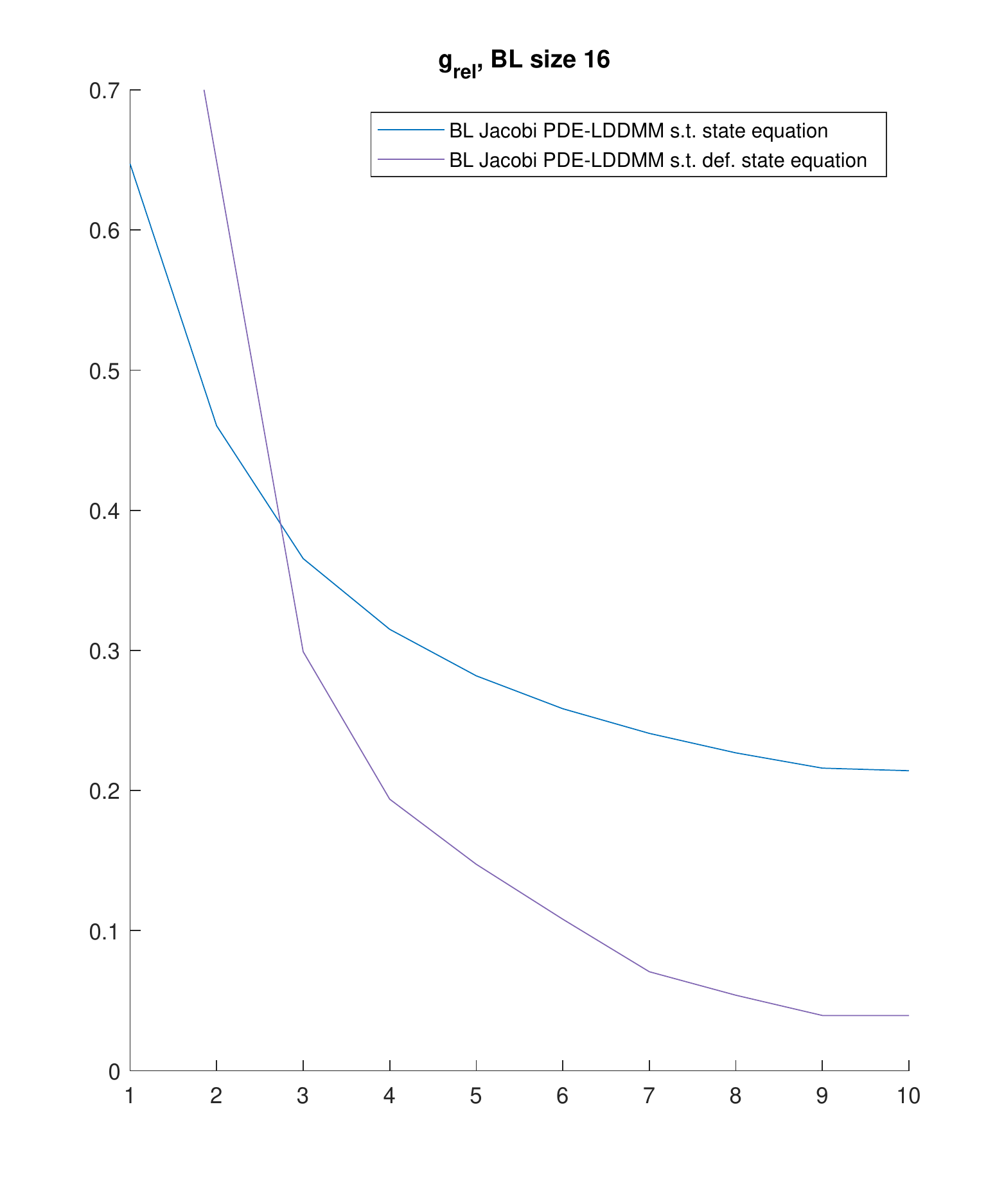} 
&
\includegraphics[width = 3 cm, height = 2.5 cm]{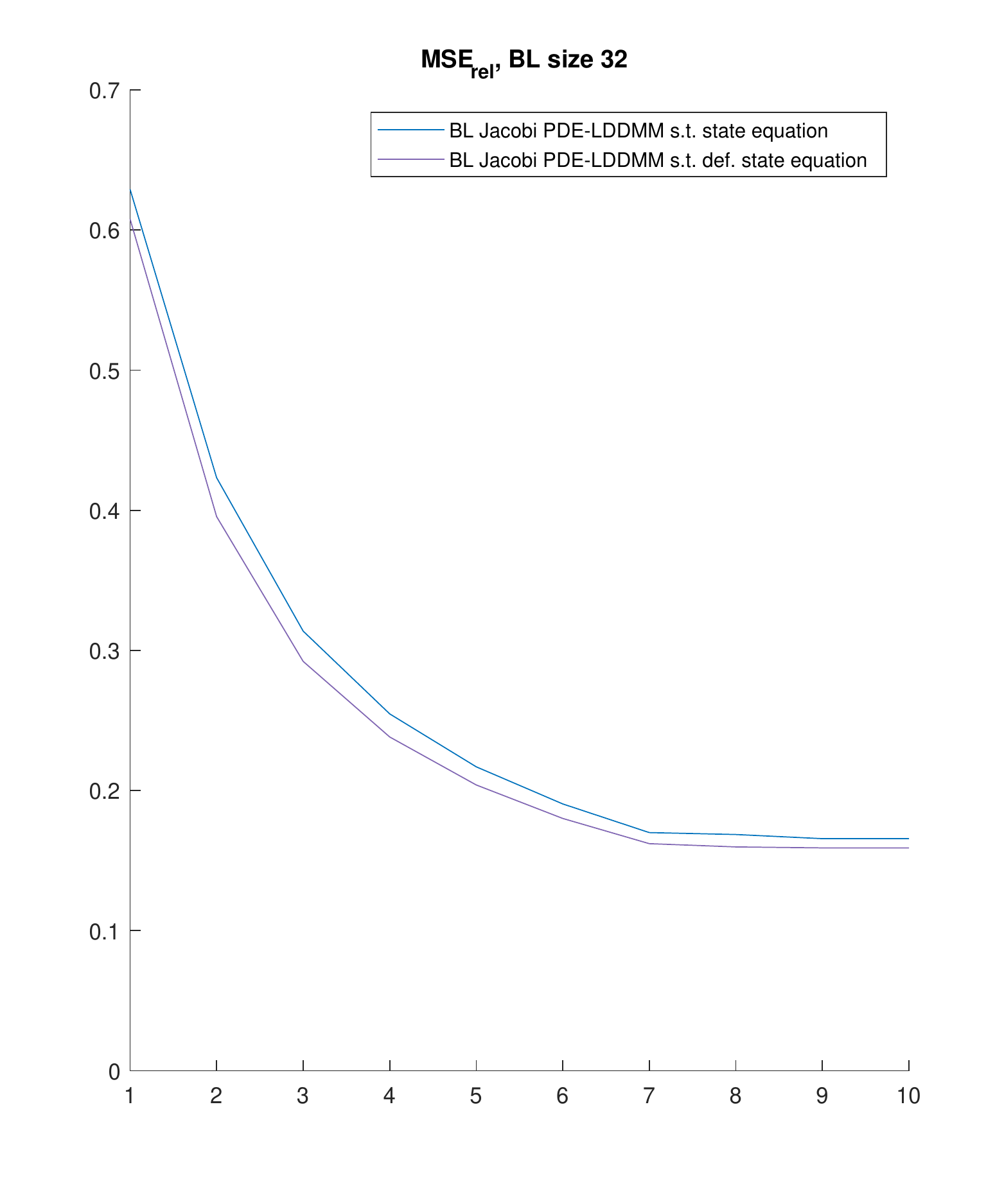} 
&
\includegraphics[width = 3 cm, height = 2.5 cm]{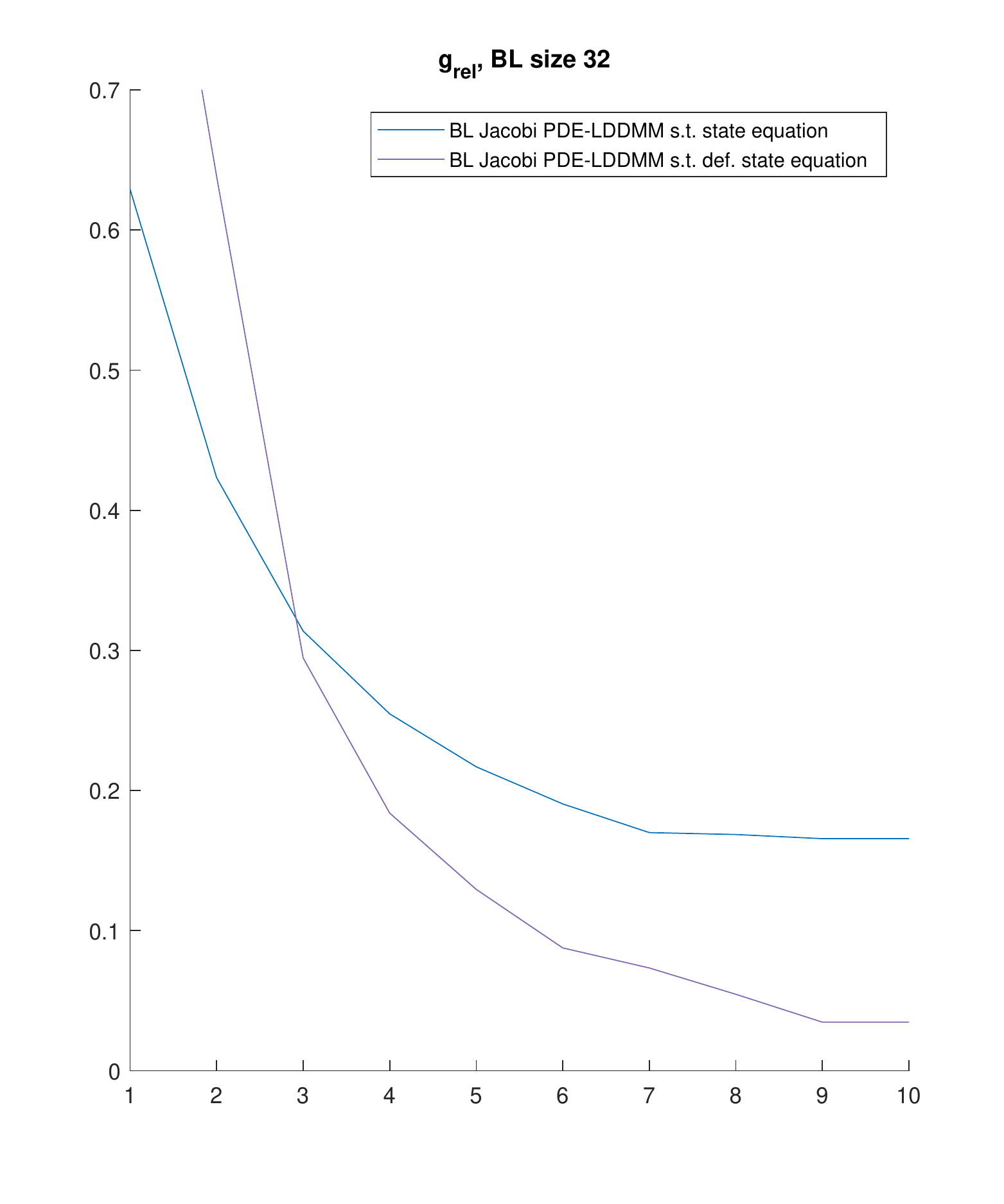} 
\\
\includegraphics[width = 3 cm, height = 2.5 cm]{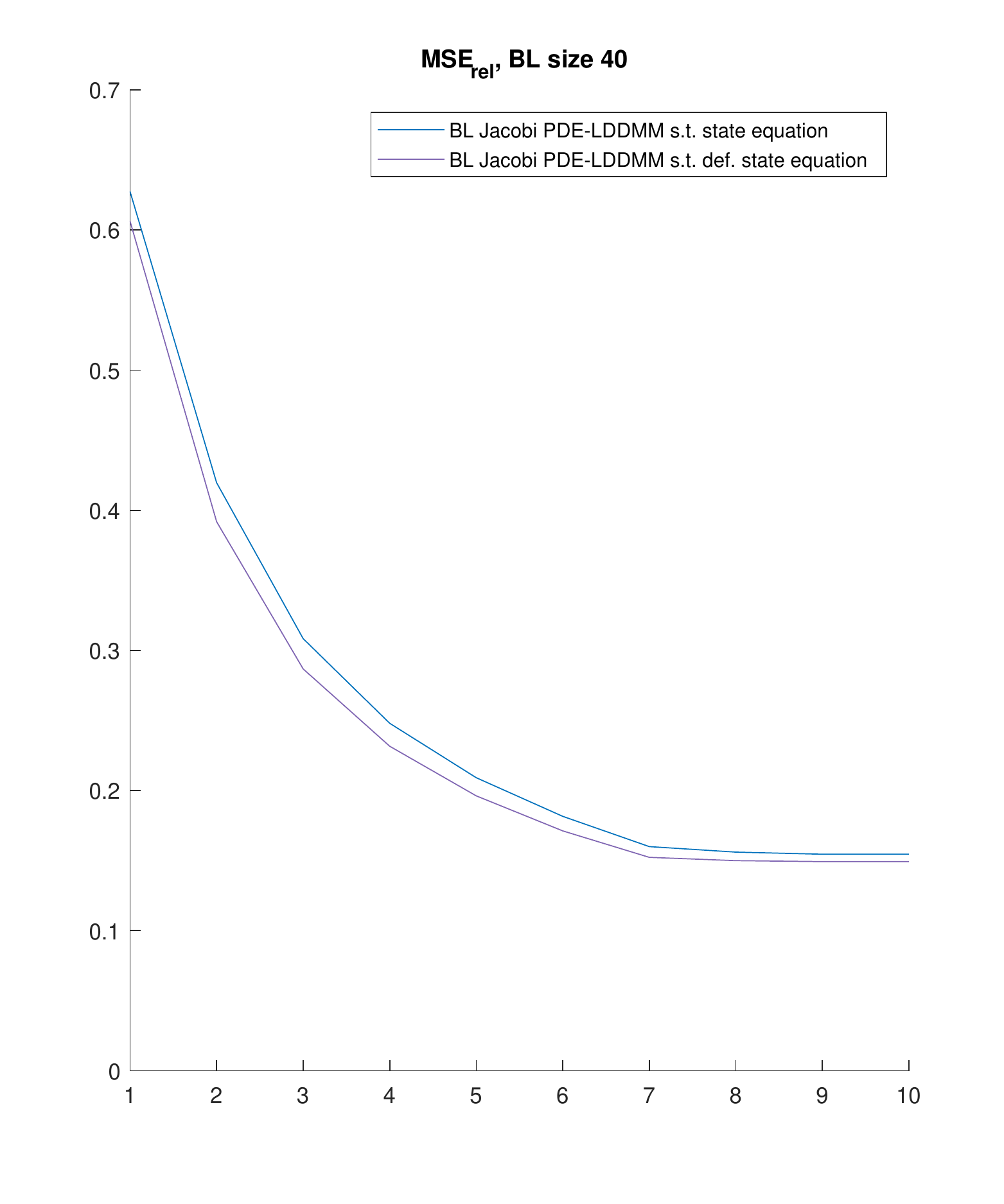} 
&
\includegraphics[width = 3 cm, height = 2.5 cm]{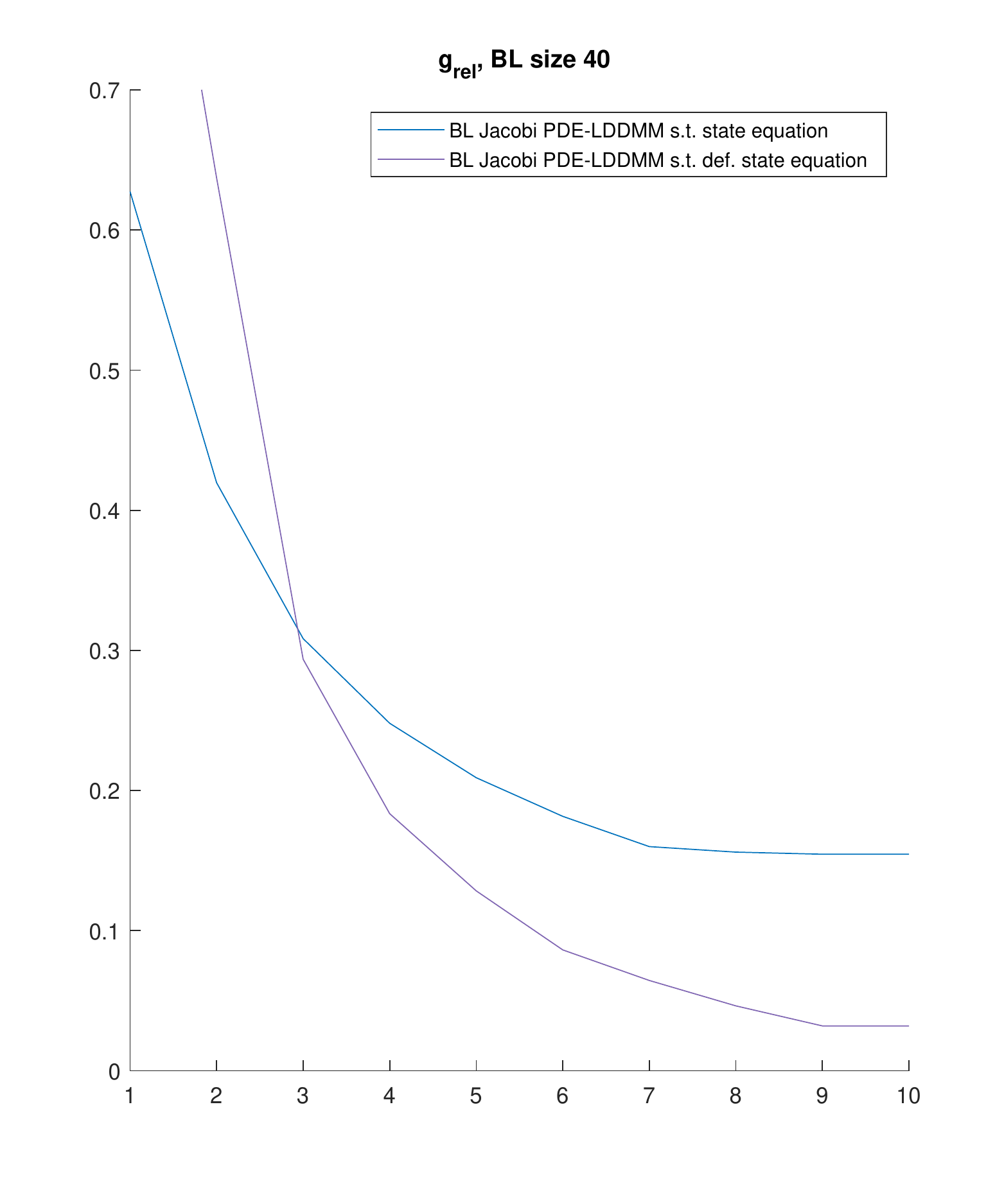} 
&
\includegraphics[width = 3 cm, height = 2.5 cm]{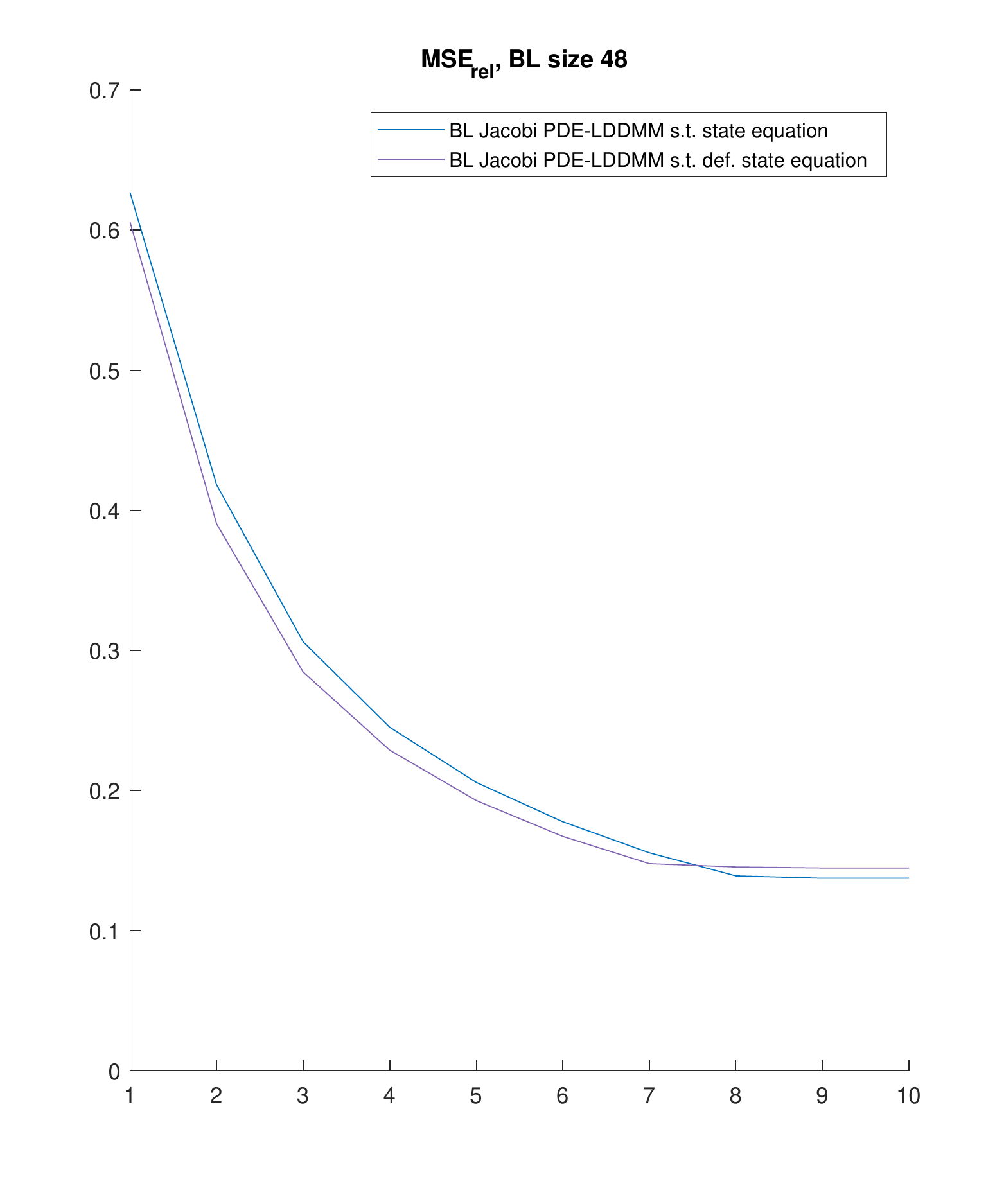} 
&
\includegraphics[width = 3 cm, height = 2.5 cm]{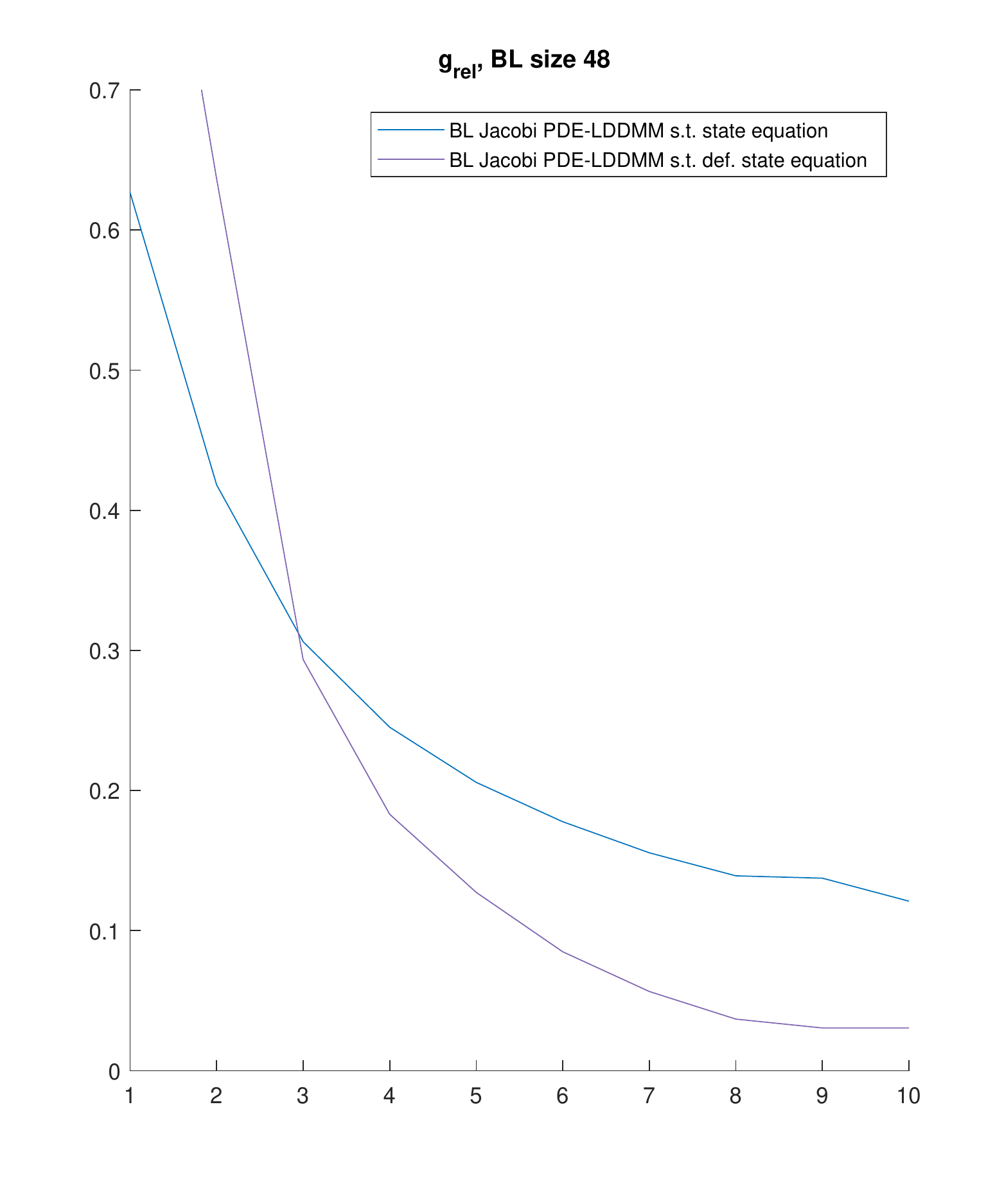} 
\\
\includegraphics[width = 3 cm, height = 2.5 cm]{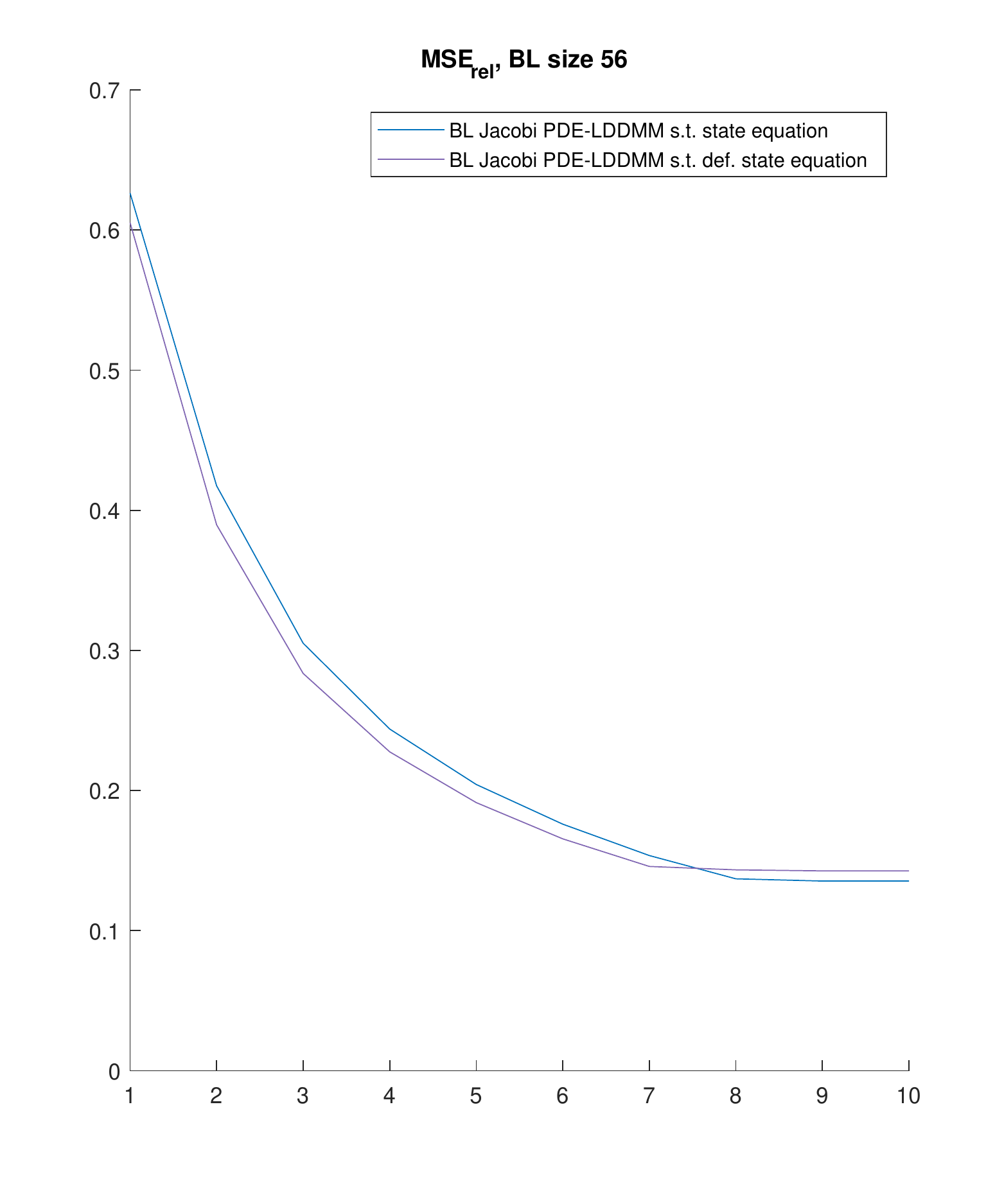} 
&
\includegraphics[width = 3 cm, height = 2.5 cm]{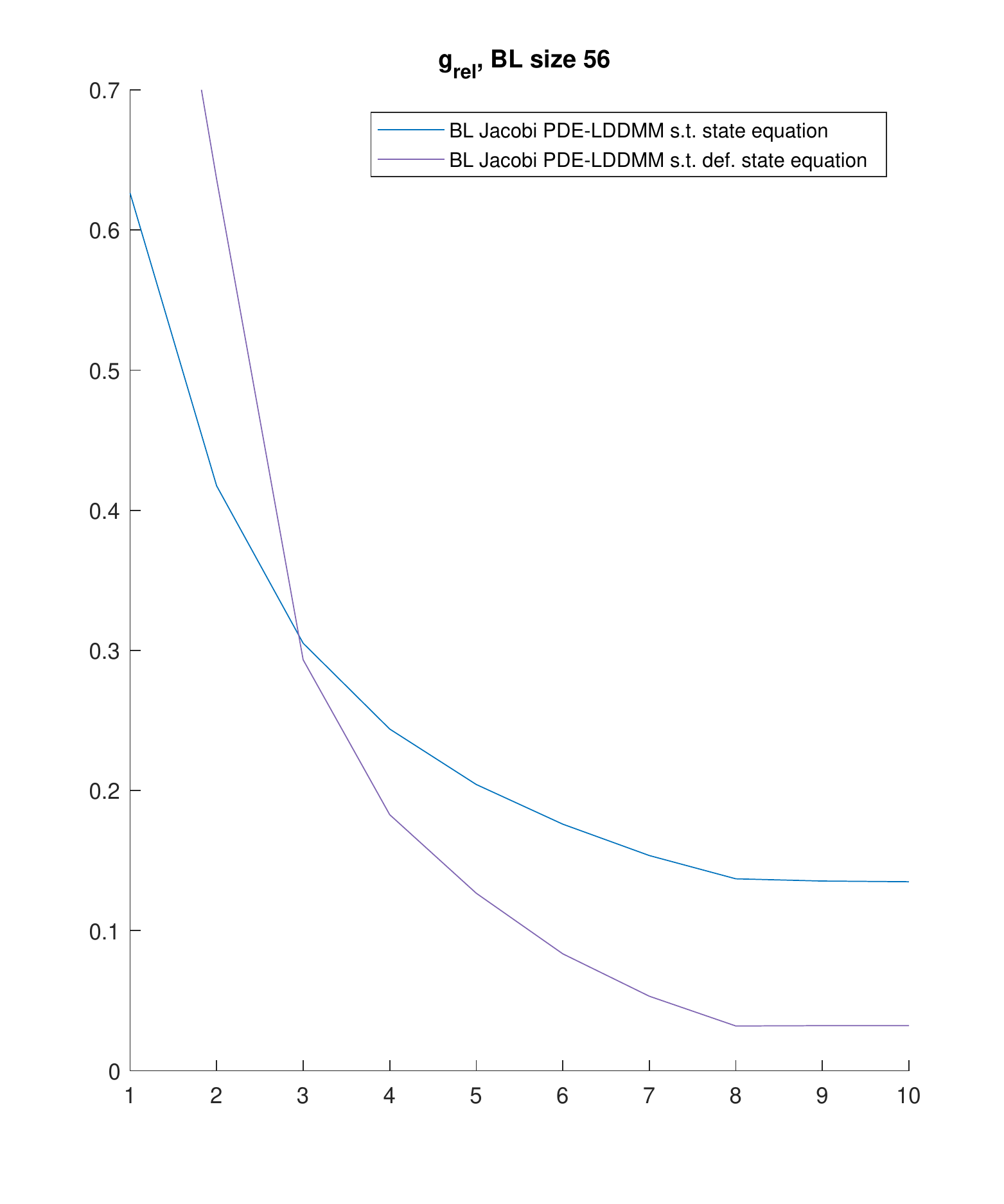} 
&
\includegraphics[width = 3 cm, height = 2.5 cm]{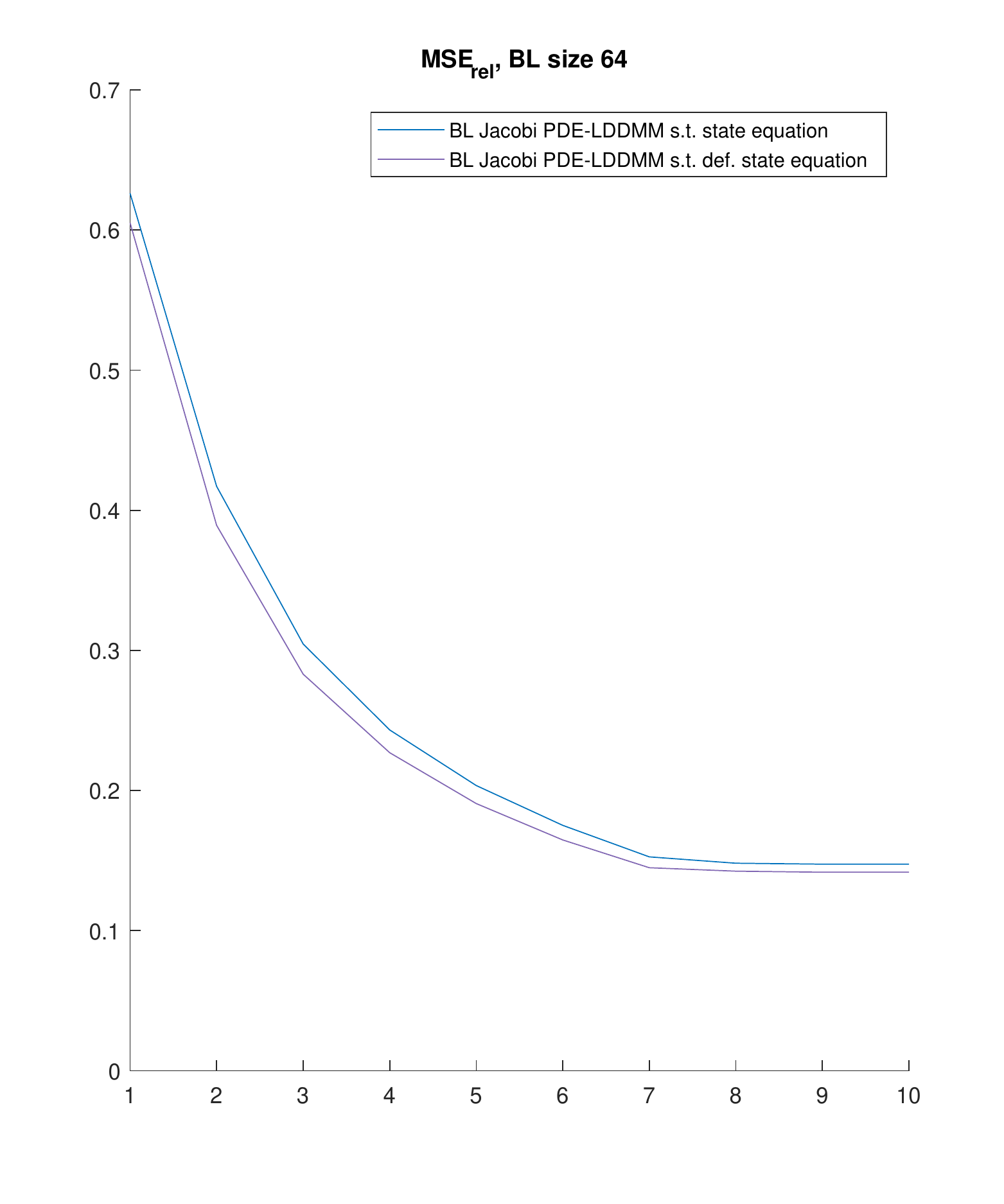} 
&
\includegraphics[width = 3 cm, height = 2.5 cm]{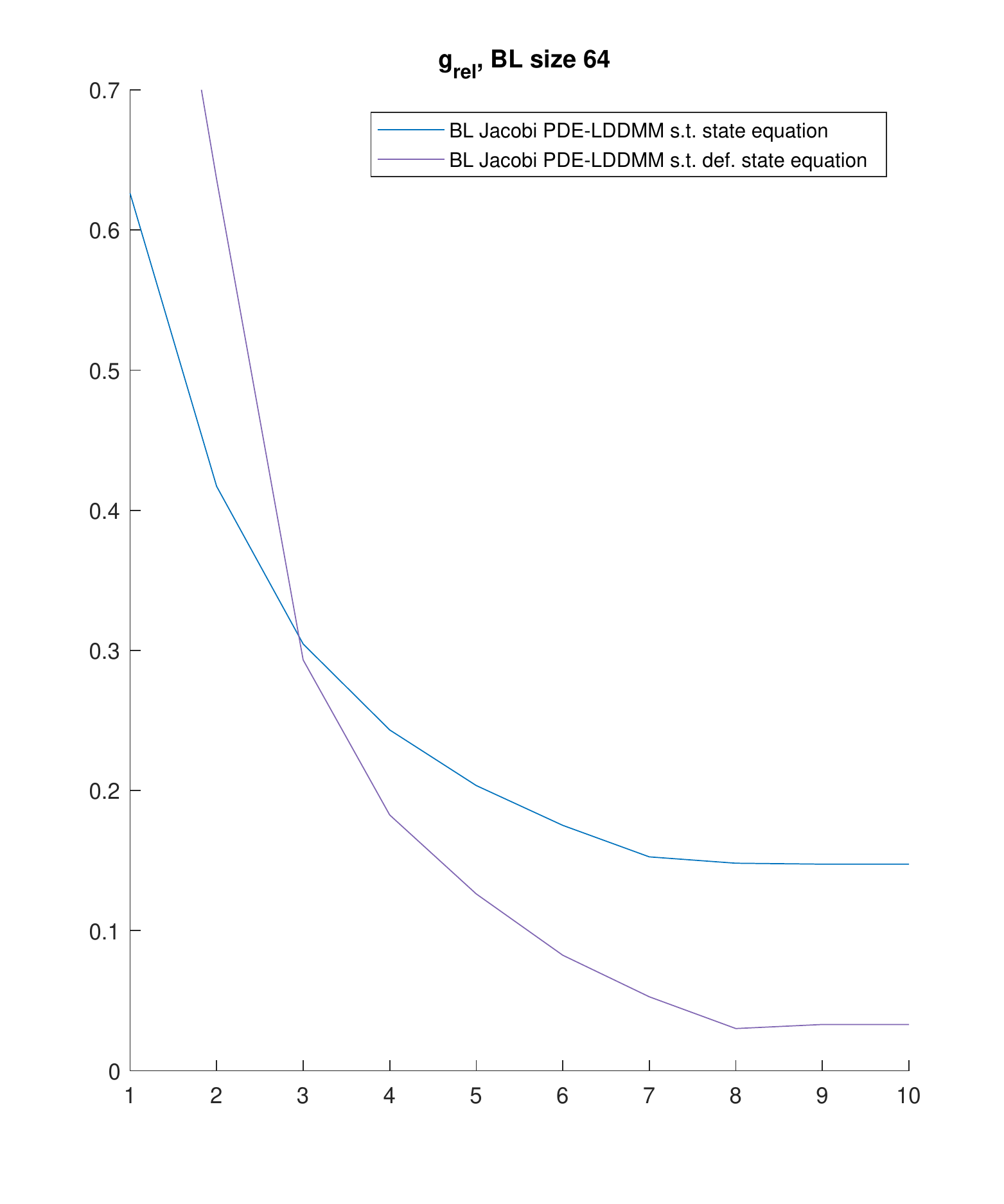} 
\\
\end{tabular}
\caption{\small Relative MSE and gradient convergence curves for the different band sizes.}
\label{fig:ConvergenceCurves}
\end{figure}

Figure~\ref{fig:ConvergenceCurves} shows the $MSE_{rel}$ and $\Vert g\Vert_{\infty,{rel}}$ convergence curves obtained during the optimization for band sizes
of 16, 32, 40, 48, 56 and 64. The figure shows that both methods converge to similar $MSE_{rel}$ values. 
However, Jacobi PDE-LDDMM subject to the deformation state equation shows smaller $\Vert g\Vert_{\infty,{rel}}$, which indicates a better convergence behavior.
Table~\ref{table:MSE} shows the numeric values after 10 iterations. 
For BL sizes of 32 both methods achieve acceptable $MSE_{rel}$ values.
Jacobi PDE-LDDMM subject to the deformation state equation slightly outperformed the method subject to the state equation.

\begin{table}[!t]
\centering

\begin{tabular}{|c||c|c|c|c|}
\hline 
 & \multicolumn{2}{|c|}{BL Jacobi PDE-LDDMM} & \multicolumn{2}{|c|}{BL Jacobi PDE-LDDMM}\\
 & \multicolumn{2}{|c|}{s.t. state eq.} & \multicolumn{2}{|c|}{s.t. def. state eq.}\\
\hline
BL size & $MSE_{rel}$ & $\Vert g\Vert_{\infty,{rel}}$ & $MSE_{rel}$ & $\Vert g\Vert_{\infty,{rel}}$ \\
\hline
16 & 21.41 & 0.03 & 21.29 & 0.03 \\
32 & 16.56 & 0.01 & 15.90 & 0.03 \\
40 & 15.45 & 0.02 & 14.92 & 0.03 \\
48 & 12.10 & 0.03 & 14.47 & 0.03 \\
56 & 13.48 & 0.03 & 14.26 & 0.03 \\
64 & 14.74 & 0.02 & 14.16 & 0.03 \\ 
\hline
\end{tabular}
\caption{\small Relative MSE and gradient achieved by the proposed methods for the different band sizes.}
\label{table:MSE}
\end{table}

\begin{figure} [!t]
\centering
\scriptsize
\begin{tabular}{ccc}
\multicolumn{3}{c}{BL Jacobi subject to state equation} \\ 
\includegraphics[width = 0.3 \textwidth]{target.pdf} 
&
\includegraphics[width = 0.3 \textwidth]{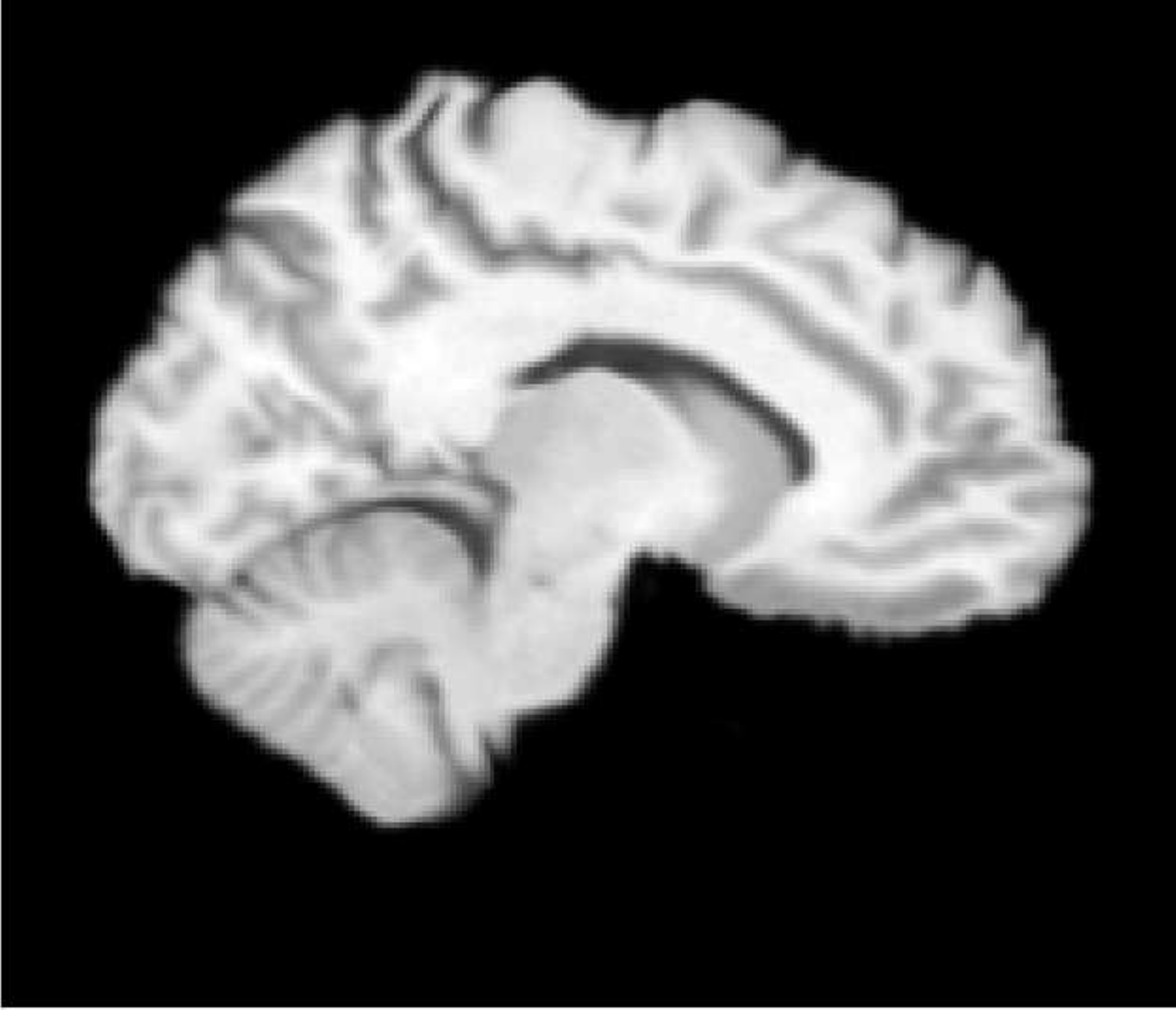} 
&
\includegraphics[width = 0.3 \textwidth]{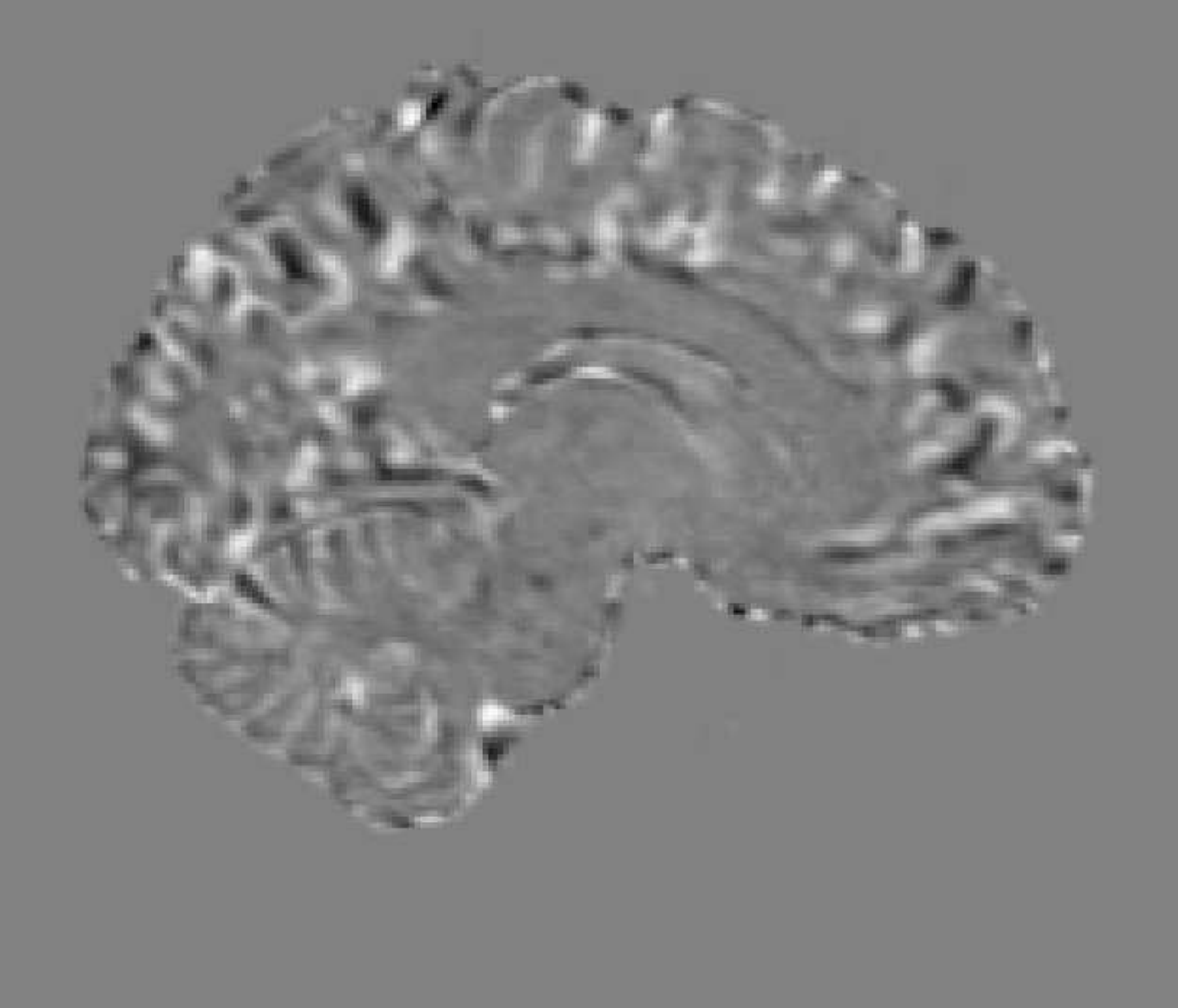} 
\\
\multicolumn{3}{c}{BL Jacobi subject to deformation state equation} \\ 
\includegraphics[width = 0.3 \textwidth]{target.pdf} 
&
\includegraphics[width = 0.3 \textwidth]{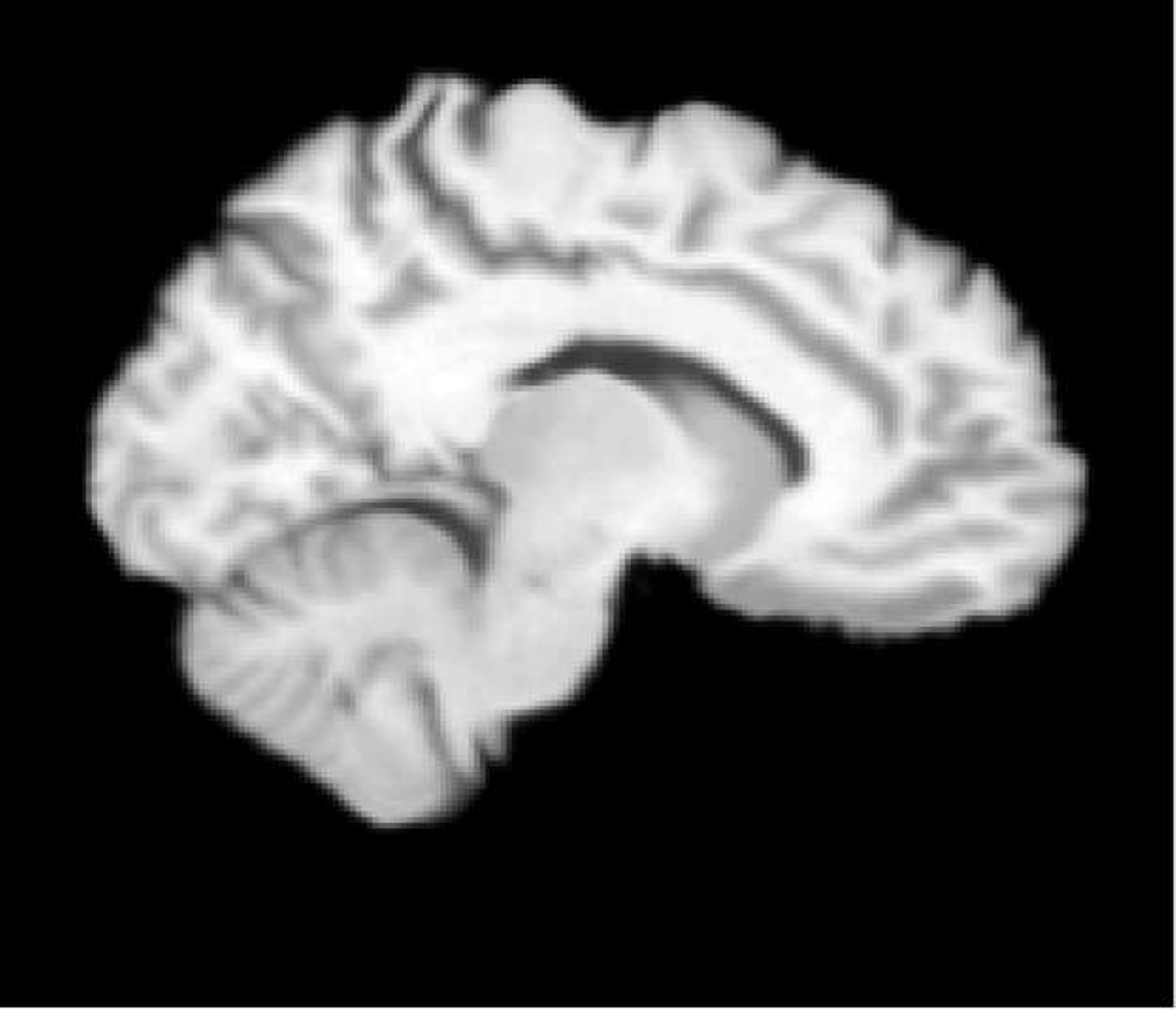} 
&
\includegraphics[width = 0.3 \textwidth]{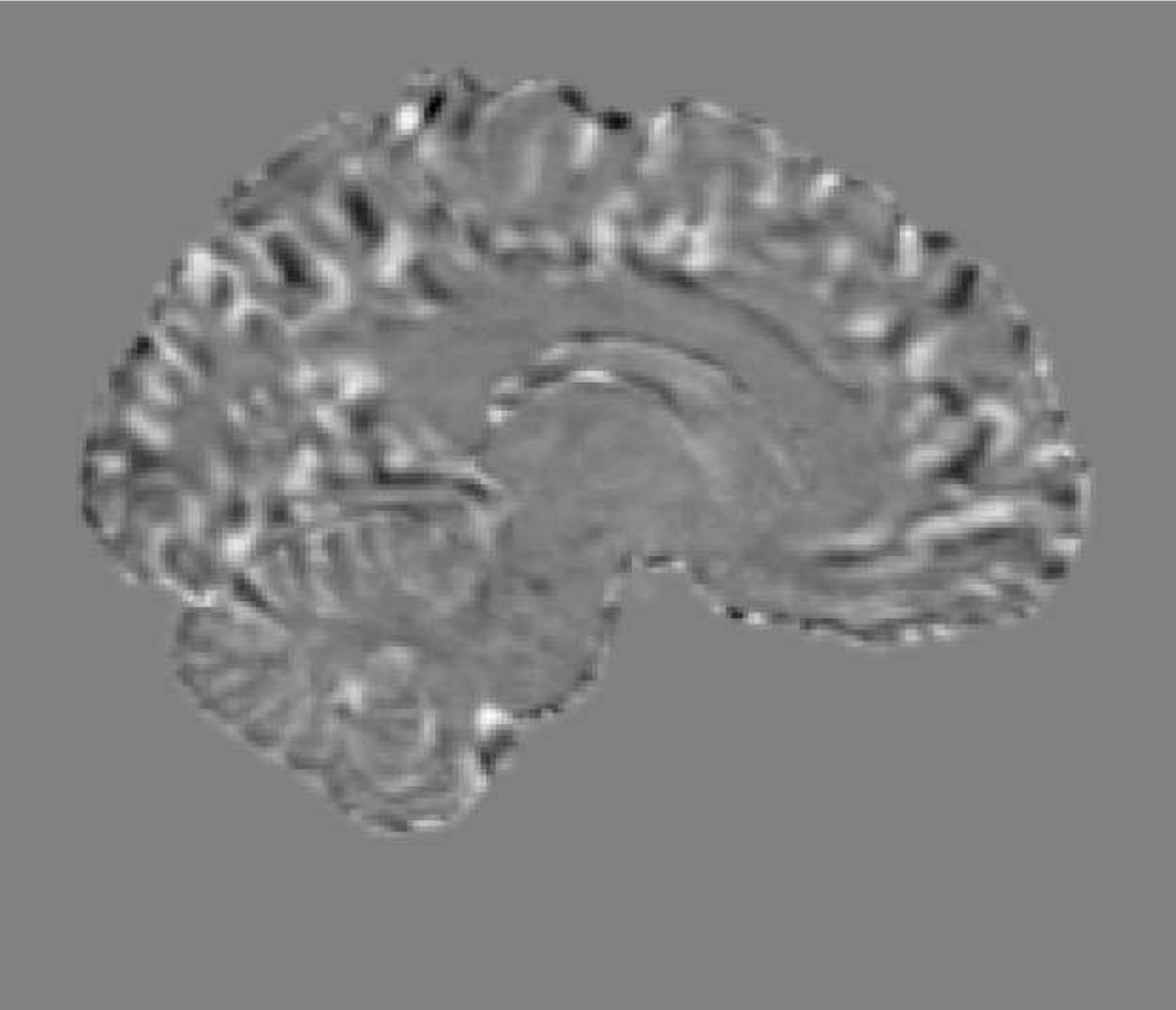} 
\end{tabular}
\caption{\small Registration results, 32x32x32. Target image, warped source, and difference after registration.}
\label{fig:QualitativeResults}
\end{figure}

Figure~\ref{fig:QualitativeResults} shows the deformed images and the differences after registration.
In the figure, it can be appreciated the accuracy achieved by the proposed methods.
The difference between the deformed images is hardly perceptible.

Finally, Table~\ref{table:ComputationalComplexity} shows the VRAM memory load and the computation time 
exhibit by the proposed methods. It should be noticed that Jacobi PDE-LDDMM subject to the state equation
in the spatial domain did not fit the memory of our available graphics card (11 GBs) and the computation
time for a downsampled example of our data was 1825.50 seconds. Therefore, the band limited vector field
parameterization definitively shows up a considerable computational saving.

\begin{table}[!t]
\centering

\begin{tabular}{|c||c|c|c|c|}
\hline 
 & \multicolumn{2}{|c|}{BL Jacobi PDE-LDDMM s.t. state eq.} & \multicolumn{2}{|c|}{BL Jacobi PDE-LDDMM s.t. def. state eq.}\\
\hline
BL size & VRAM (MBS) & total time (s) & VRAM (MBS) & total time (s) \\
\hline
16 & 3319 & 892.56 & 1745 & 694.43 \\ 

32 & 3439 & 894.38 & 1927 & 701.14 \\ 

40 & 3579 & 845.57 & 2137 & 653.20 \\ 

48 & 3803 & 869.67 & 2469 & 693.83 \\ 

56 & 4097 & 1154.58 & 2907 & 1022.55 \\ 

64 & 4403 & 1576.76 & 3393 & 1492.98 \\ 
\hline
\end{tabular}
\caption{\small Computational complexity of the proposed methods. }
\label{table:ComputationalComplexity}
\end{table}

\section*{Acknowledgements}
This work was partially supported by Spanish research grant TIN2016-80347-R.

\bibliographystyle{splncs}
\bibliography{abbsmall.bib,Diffeo.bib,OpticalFlow.bib}

\begin{thebibliography}{10}

\bibitem{Sotiras_13}
Sotiras, A., Davatzikos, C., Paragios, N.:
\newblock Deformable medical image registration: A survey.
\newblock IEEE Trans. Med. Imaging \textbf{32(7)} (2013)  1153 -- 1190

\bibitem{Miller_04}
Miller, M.I.:
\newblock Computational anatomy: shape, growth, and atrophy comparison via
  diffeomorphisms.
\newblock Neuroimage \textbf{23} (2004)  19--33

\bibitem{Hart_09}
Hart, G.L., Zach, C., Niethammer, M.:
\newblock An optimal control approach for deformable registration.
\newblock Proc. of the IEEE Conference on Computer Vision and Pattern
  Recognition (CVPR'09) (2009)

\bibitem{Beg_05}
Beg, M.F., Miller, M.I., Trouve, A., Younes, L.:
\newblock Computing large deformation metric mappings via geodesic flows of
  diffeomorphisms.
\newblock Int. J. Comput. Vision \textbf{61 (2)} (2005)  139--157

\bibitem{Vialard_12}
Vialard, F.X., Risser, L., Rueckert, D., Holm, D.D.:
\newblock Diffeomorphic atlas estimation using geodesic shooting on volumetric
  images.
\newblock Annals of the BMVA \textbf{2012} (2012)  1 -- 12

\bibitem{Mang_15}
Mang, A., Biros, G.:
\newblock An inexact {N}ewton-{K}rylov algorithm for constrained diffeomorphic
  image registration.
\newblock SIAM J. Imaging Sciences \textbf{8(2)} (2015)  1030--1069

\bibitem{Hernandez_18_ArxivMICCAI}
Hernandez, M.:
\newblock {PDE}-constrained {LDDMM} via geodesic shooting and inexact
  {G}auss-{N}ewton-{K}rylov optimization using the incremental adjoint {J}acobi
  equations.
\newblock ArXiv:1807.04638 (2018)

\bibitem{Hernandez_18}
Hernandez, M.:
\newblock Band-{L}imited {S}tokes {L}arge {D}eformation {D}iffeomorphic
  {M}etric {M}apping.
\newblock IEEE J. of Biom. and Health Inf. (2018)

\bibitem{Hernandez_18_ArxivECCV}
Hernandez, M.:
\newblock {N}ewton-{K}rylov {PDE}-constrained {LDDMM} in the space of
  band-limited vector fields.
\newblock ArXiv:1807.05117 (2018)

\bibitem{Holm_98}
Holm, D.D., Marsden, J.E., Ratiu, T.S.:
\newblock The {E}uler-{P}oincare equations and semidirect products with
  applications to continuum theories.
\newblock Adv. in Math. \textbf{137} (1998)  1 -- 81

\bibitem{Younes_07}
Younes, L.:
\newblock Jacobi fields in groups of diffeomorphisms and applications.
\newblock Q. Appl. Math. \textbf{65} (2007)  113 -- 134

\bibitem{Zhang_15}
Zhang, M., Fletcher, P.T.:
\newblock Finite-dimensional {L}ie algebras for fast diffeomorphic image
  registration.
\newblock Proc. of International Conference on Information Processing and
  Medical Imaging (IPMI'15), Lecture Notes in Computer Science (2015)

\bibitem{Bullo_95}
Bullo, F.:
\newblock Invariant affine connections and controllability on {L}ie groups.
\newblock Techical report for Geometric Mechanics, California Institute of
  Technology (1995)

\bibitem{Polzin_16}
Polzin, T., Niethammer, M., Heinrich, M.P., Handels, H., Modersitzki, J.:
\newblock Memory efficient {LDDMM} for lung {CT}.
\newblock (2014)  28--36

\end{thebibliography}

\end{document}